\PassOptionsToPackage{reqno}{amsmath}
\documentclass[final,hidelinks,onefignum,onetabnum]{siamart251216}

\usepackage{amssymb,mathtools}
\usepackage[ntheorem]{empheq}
\usepackage{enumitem}
\usepackage{xcolor}
\usepackage{subcaption}
\usepackage{float}
\usepackage{graphicx}
\usepackage{algorithm}
\usepackage{algpseudocode}
\usepackage{booktabs,tabularx}
\usepackage{placeins}
\usepackage{cite}
\definecolor{ink}{RGB}{24,48,75}
\definecolor{sea}{RGB}{0,112,128}
\definecolor{pale}{RGB}{232,244,246}
\definecolor{warm}{RGB}{179,92,36} 

\theoremstyle{plain}
\theoremheaderfont{\normalfont\scshape}
\theorembodyfont{\normalfont\itshape}
\theoremseparator{.}
\theoremsymbol{}

\renewtheorem{lemma}{Lemma}
\newtheorem{assumption}{Assumption}

\renewtheorem{corollary}{Corollary}[section]
\theorembodyfont{\normalfont}
\renewtheorem{definition}[theorem]{Definition}
\theoremheaderfont{\normalfont\itshape}
\newtheorem{remark}{Remark}

\newcommand{\R}{\mathbb{R}}
\newcommand{\norm}[1]{\left\|#1\right\|}

\providecommand{\alert}[1]{#1}

\def\call{\mathcal{L}}

\DeclareMathOperator*{\argmin}{arg\,min}
\newcommand{\NN}{\mathcal N_\theta}
\newcommand{\NNa}{\mathcal N_\theta^\circ}

\newcommand{\PB}{\mathrm{PBU}}

\numberwithin{equation}{section}
\numberwithin{theorem}{section}
\numberwithin{lemma}{section}
\numberwithin{assumption}{section}
\numberwithin{remark}{section}

\definecolor{proofred}{RGB}{190,0,0}

\usepackage{tikz}
\usetikzlibrary{arrows.meta,positioning,calc}

\headers{Stable admissible neural classes}{L. Ling, S. Ma, and Y. Sun}
\title{\texorpdfstring{Stable admissible neural classes and convergent physics-informed neural networks}{Stable admissible neural classes and convergent physics-informed neural networks}\thanks{\textbf{Funding.}
This work was supported by the General Research Fund
(GRF No. 12301824, 12301526) of the Hong Kong Research Grants Council,
the Guangdong and Hong Kong Universities ``1+1+1'' Joint Research
Collaboration Scheme (project no. 2025A0505000014),
the National Natural Science Foundation of China
(Young Scientists Fund, Category C, Grant No. 12601787)
and the Guangdong Basic and Applied Basic Research Foundation
(Grant No. 2026A1515012472).}}

\author{Leevan Ling\thanks{Department of Mathematics, Hong Kong Baptist University,
Hong Kong (\email{lling@hkbu.edu.hk}, \href{mailto:shu_ma@hkbu.edu.hk}{shu\_ma@hkbu.edu.hk},
\email{yzsun95@hkbu.edu.hk}).}
\and Shu Ma\footnotemark[2]\kern0.35em\thanks{Corresponding author.}
\and Yizhong Sun\footnotemark[2]}

\begin{document}

\maketitle

\begin{abstract}
Residual-based methods, including physics-informed neural networks (PINNs), physics-informed Kolmogorov--Arnold networks, and some physics-informed neural operator architectures, approximate PDE solutions by minimizing a discrete least-squares loss formed from the PDE and boundary residuals over a prescribed trial class. This structure is shared by classical least-squares Kansa collocation, whose kernel trial spaces admit an established convergence theory. For neural trial classes, approximation power alone does not ensure that discrete residuals control the continuous solution error. For linear elliptic Dirichlet problems, we extend the kernel-based collocation framework to neural trial classes through elliptic regularity, discrete stability, and consistency. Elliptic regularity supplies the PDE estimate relating solution error to continuous residuals. Discrete stability identifies admissible trial subclasses in which continuous residuals are controlled by a computable discrete loss. This step also quantifies the required number of collocation points. Furthermore, the stability analysis naturally demands a scale-dependent, physics-balanced boundary weight. For the constrained problem, consistency requires accurate comparison functions within the same admissible subclass, a requirement not supplied by neural approximation estimates alone. Specializing to multilayer perceptrons with tanh activations (tanh-MLPs), we construct admissible comparison networks and establish convergence for the constrained residual minimization problem. We further obtain conditional convergence of unconstrained physics-balanced PINNs based on tanh-MLPs under joint network and sampling refinement.
\end{abstract}


% \begin{MSCcodes}
% 65N12, 65N35, 65D12, 68T07
% \end{MSCcodes}

\section{\texorpdfstring{Motivation and analytical framework}{Motivation and analytical framework}}
Residual-based methods approximate solutions of partial
differential equations by minimizing equation and boundary residuals over
a prescribed trial class $\NN$ with an assigned design scale $h_\theta$.
Its BII-admissible subclass is denoted by $\NNa$.
For the elliptic problem considered here, regularity estimates bound the solution error in terms of continuous residual norms. We then use kernel-based collocation to motivate the construction of stable admissible neural classes.
\subsection{\texorpdfstring{Elliptic regularity and residual minimization}{Elliptic regularity and residual minimization}}
\label{sec:continuous-residual-framework}
Consider the linear elliptic Dirichlet problem on $\Omega\subset\R^d$ ($d\le3$):
\begin{equation}\label{eq:PDE}
\call u=f \quad \text{in }\Omega,
\qquad u=g \quad \text{on }\partial\Omega,
\end{equation}
where
\begin{equation}\label{eq:L}
\call u
:=
-\sum_{i,j=1}^d
\frac{\partial}{\partial x_j}
\left(a^{ij}(x)\frac{\partial u}{\partial x_i}\right)
+\sum_{i=1}^d b^i(x)\frac{\partial u}{\partial x_i}
+c(x)u.
\end{equation}
We make the following assumptions.
\begin{assumption}\label{ass:pde}
Let $m\ge2$.
\begin{enumerate}[label=(\arabic*),ref=\theassumption(\arabic*)]
\item\label{ass:Omega}
We assume that the bounded domain
$\Omega$ has a piecewise $C^m$--boundary $\partial\Omega$ so that $\Omega$
is a Lipschitz domain and satisfies an interior cone condition.
\item
We assume that the functions $f$ and $g$ are smooth enough to admit
a unique \emph{classical} solution $u^*\in H^m(\Omega)$ of \eqref{eq:PDE}.
\item\label{ass:L}
The operator $\call$ in \eqref{eq:L} is strongly elliptic, with
coefficients in $W^{m,\infty}(\Omega)$. 
\end{enumerate}
\end{assumption}

The trace theorem \cite{Wloka1987Partial} can be applied and
we can define the trace operator
$\operatorname{Tr}:H^s(\Omega)\longrightarrow H^{s-1/2}(\partial\Omega)$ for $\tfrac12<s\le m$, such that $\operatorname{Tr}u=u|_{\partial\Omega}$ for all $u\in C^m(\bar\Omega)$,
with a continuous right-inverse linear extension operator on the
corresponding trace space.
The smoothness assumption also allows boundary coordinate patches with a
subordinate partition of unity. Each patch can be mapped to the unit ball
in $\R^{d-1}$ by a $C^m$--diffeomorphism.
Under these assumptions, $\call:H^m(\Omega)\to H^{m-2}(\Omega)$
is bounded \cite{GieslWendland2007Meshless}, with
\begin{equation}\label{eq:Lu<uO}
\|\call u\|_{H^{m-k-2}(\Omega)}\leq C_{\Omega,\call} \|u\|_{H^{m-k}(\Omega)},
\;k\in\mathbb{N}_0,\quad k\le m-2,
\end{equation}
for all $u\in H^m(\Omega)$. Moreover, the following boundary
regularity estimate \cite{Jost2007Partial} holds
\begin{align}
\| u\|_{H^{k+2}(\Omega)}
&\le C_{\Omega,\call,k} \left( \| \call u \|_{H^k(\Omega)}
+ \| u \|_{H^{k+3/2}( \partial\Omega)} 
\right),\; 0\leq k \leq m-2,
\label{eq:bdyreg}
\end{align}
for all $u\in H^m(\Omega)$ with $C_{\Omega,\call,k}$ depending on $\Omega$, the ellipticity constant of $\call$, and $k\geq 0$.

\medskip
The preceding regularity estimate bounds the solution error in terms of continuous residual norms. In strong-form collocation, we evaluate these residuals at finitely many interior and boundary points.
This formulation includes multilayer perceptron (MLP)-based
physics-informed neural networks (PINNs)~\cite{RaissiPerdikarisKarniadakis2019,
SanjalaweEtAl2026PINNReview}, physics-informed Kolmogorov--Arnold networks
(PI-KANs)~\cite{Noorizadegan+WangETAL-MakiGausnetwreli:26}, and some physics-informed neural operator architectures.
In these methods, $\NN$ is a neural network class
consisting of functions represented by a prescribed architecture as its
parameters vary.  Such a class is generally not a linear space.
Shallow networks with fixed hidden parameters yield linear trial
spaces~\cite{MaoXuXu2026LinearizedPDEs}.
The central issue is whether minimizing sampled residuals over $\NN$
yields accurate approximations to the PDE solution.

Neural approximation results provide accurate approximants to PDE
solutions~\cite{GuhringRaslan2021,DeRyckLanthalerMishra2021,
MaoSiegelXu2026Radon} and solution operators~\cite{CaoGuoHuangLiu2026OperatorBarron}.
To control the solution error, existing analyses use continuous stability
estimates in terms of residual norms~\cite{MishraMolinaro2023GeneralizationPINN,
ZeinhoferMasriMardal2025UnifiedPINN}.
Sampling, quadrature, or generalization estimates then relate these residual
norms to the discrete objective~\cite{ShinDarbonKarniadakis2020ConvergencePINN,
ShinZhangKarniadakis2023ResidualMinimization}.
For a general neural network class, however, approximation power alone
does not ensure uniform control of continuous residual norms by their
sampled counterparts. 
Convergence therefore requires both sufficiently accurate comparison functions within
the trial class and discrete stability estimates that control the solution
error by sampled residuals.
We next review how stability and consistency are combined in
kernel-based collocation~\cite{CheungLingSchaback2018H2,
ChenLingYun2025KernelStability}.

\subsection{\texorpdfstring{Review of stability and consistency in kernel-based collocation}{Review of stability and consistency in kernel-based collocation}}

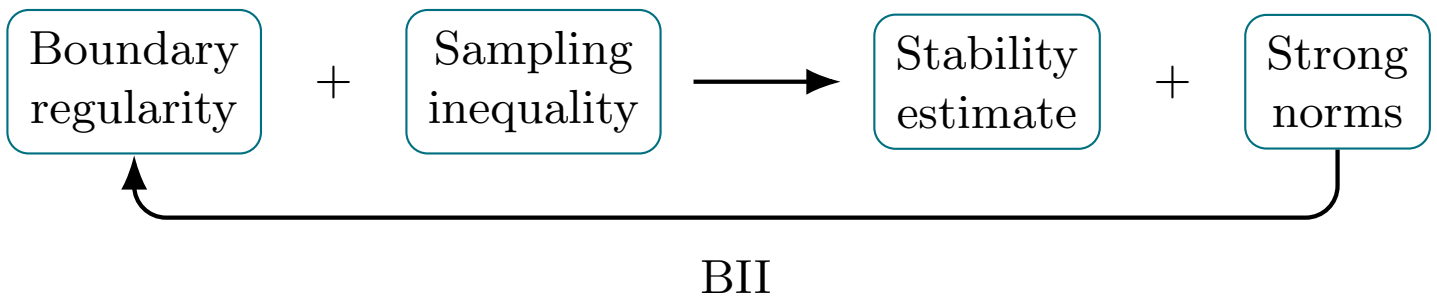
\begin{figure}[t]
\centering
\begin{tikzpicture}[node distance=0.20cm,>=Latex,
every node/.style={draw=sea,rounded corners,align=center,
minimum height=0.75cm,inner sep=4pt,font=\small}]
\node (a) {Boundary\\regularity};
\node[right=of a,draw=none] (plus1) {$+$};
\node[right=of plus1] (b) {Sampling\\inequality};
\node[right=1.4cm of b] (c) {Stability\\estimate};
\node[right=of c,draw=none] (plus2) {$+$};
\node[right=of plus2] (d) {Strong\\norms};
\draw[->,thick,shorten <=6pt,shorten >=6pt] (b)--(c);
\draw[->,thick,rounded corners=6pt]
(d.south) -- ++(0,-0.45) -|
node[pos=0.25,below,draw=none,font=\small] {BII} (a.south);
\end{tikzpicture}
\vspace{-10pt}
\caption{Discrete stability in least-squares kernel-based collocation. Regularity and sampling inequalities bound the $H^2$ norm of trial functions by sampled residual terms and strong norm remainders. The return arrow shows how a Bernstein-type inverse inequality (BII) enables absorption under sufficiently dense sampling.}
\vspace{-10pt}
\label{fig:bii-discrete-stability}
\end{figure}

\begin{figure}[t]
\begingroup\color{black}
\parindent=0pt
\centering
\begin{tikzpicture}[x=.98\linewidth,y=1cm,
every node/.style={text=black,inner sep=0pt,anchor=north west,align=left},
ingredient/.style={draw=sea,line width=.65pt,rounded corners=3pt},
heading/.style={font=\footnotesize\bfseries,anchor=north,align=center},
description/.style={font=\footnotesize,anchor=west,align=justify}]
% \foreach \x in {0,.34,.68} {\draw[ingredient] (\x,0) rectangle +(0.32,2.30);}
\foreach \x in {0} {\draw[ingredient] (\x,0) rectangle +(0.3,2.30);}
\foreach \x in {.315} {\draw[ingredient] (\x,0) rectangle +(0.33,2.30);}
\foreach \x in {.66} {\draw[ingredient] (\x,0) rectangle +(0.34,2.30);}
\node[heading] at (.15,2.15) {Regularity};
\node[heading] at (.48,2.15) {Stability};
\node[heading] at (.83,2.15) {Consistency};
\node[description,text width=.265\linewidth] at (.015,.98)
{PDE regularity controls solution norms through the data $f$ and $g$.};
\node[description,text width=.294\linewidth] at (.33,.98)
{Sampling and inverse inequalities bound
continuous error norms by discrete loss.};
\node[description,text width=.304\linewidth] at (.675,.98)
{A good comparison function from the same
trial space yield convergence rates
when combined with stability.};
\end{tikzpicture}
\par\endgroup
\caption{Three ingredients for convergence.}
\label{fig:three-ingredients-convergence}
\vspace{-15pt}
\end{figure}
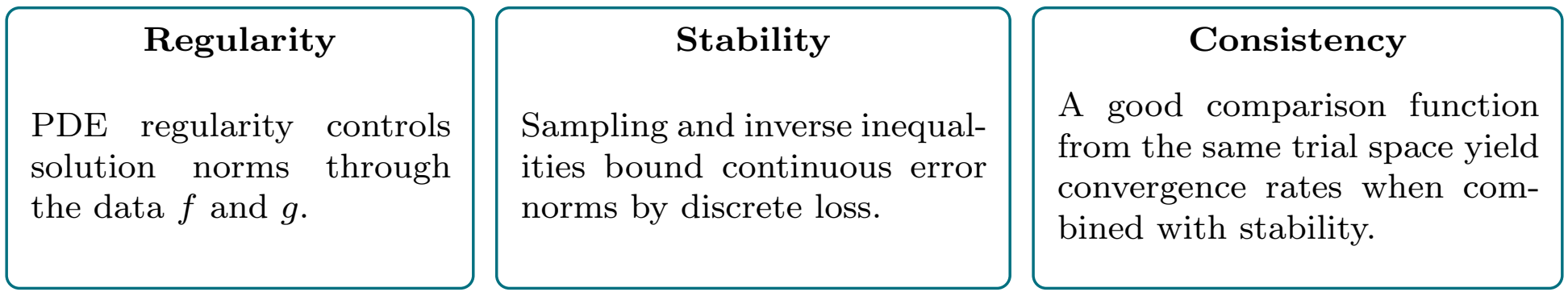

Kernel-based collocation~\cite{Kansa1990PartII} provides one way to combine
accurate approximation with discrete stability within the same trial space.
For a finite center set $Z\subset\Omega$, the kernel trial space is the
finite-dimensional reproducing kernel Hilbert space (RKHS)
\begin{align}
V_Z:=\operatorname{span}\{\Phi(\cdot,z):z\in Z\}
\subset\mathcal H_\Phi(\Omega)=H^m(\Omega),\qquad m>d/2,
\end{align}
with a translation-invariant positive definite kernel
$\Phi(x,z)=\phi(x-z)$ with
$\widehat\phi(\omega)\simeq(1+|\omega|^2)^{-m}$
and inherited inner product~\cite{wendland2004scattered,CheungLingSchaback2018H2}. Let $X\subset\Omega$ and $Y\subset\partial\Omega$ be interior and boundary
collocation points. Set $\|\phi\|_{\ell^2(X)}^2 :=\sum_{x\in X}|\phi(x)|^2$.
Write $h_Z,h_X,h_Y$ for the fill distances of $Z,X,Y$ in
$\Omega,\Omega,\partial\Omega$, respectively, and $q_Z$ for the separation
radius of $Z$. For quasi-uniform centers, we have
$q_Z\simeq h_Z$~\cite{wendland2004scattered}.

Let $u\in V_Z$, with $m>2+d/2$.
Applying PDE regularity \cite{Jost2007Partial} and sampling \cite{ArcangeliLopezTorrens2012Sampling}:
\begin{align}
&\|u\|_{H^2(\Omega)}^2
\lesssim \|\mathcal Lu\|_{L^2(\Omega)}^2
+\|u\|_{H^{3/2}(\partial\Omega)}^2\notag
\\[0.8em]
&\phantom{\|u\|_{H^2(\Omega)}^2}
\lesssim  
\underbrace{
\alert{h_X^{d}\|\mathcal Lu\|_{\ell^2(X)}^2
+h_Y^{d-4}\|u\|_{\ell^2(Y)}}^2
}_{\displaystyle  \text{Weighted SOS residual }  J(u)}
+\underbrace{\bigl(h_X^{2m-4}+h_Y^{2m-4}\bigr)
\|u\|_{H^m(\Omega)}^2}_{\displaystyle\text{strong-norm remainder}}. 
\end{align}
For an $H^m(\Omega)$-Sobolev kernel, with $m>2+d/2$,
the trial space $V_Z$ satisfies a Bernstein inverse
inequality  (BII) \cite{ArcangeliLopezTorrens2012Sampling,NarcowichWardWendland2006Bernstein}
\begin{align}
\|u\|_{H^m(\Omega)}\le C' q_Z^{\,2-m}\|u\|_{H^2(\Omega)}     \qquad
\text{for every }u\in V_Z,
\end{align}
Under suitable oversampling condition \cite{CheungLingSchaback2018H2}, the remainder is absorbed and gives
% $C' \bigl(h_X^{2m-4}+h_Y^{2m-4}\big) q_Z^{\,4-2m} \le \tfrac12$ 
\begin{flalign}\label{kernel-stability}
\makebox[12em][l]{(Stability)}
&\|u\|_{H^2(\Omega)}^2\lesssim J(u)
\qquad\text{for every }u\in V_Z,&&
\end{flalign}
which yields discrete stability on $V_Z$, see Figure~\ref{fig:bii-discrete-stability}.
The method of comparison now combines this discrete stability estimate
with accurate comparison functions from the same trial space $V_Z$.
For a least-squares minimizer $u_Z\in\operatorname*{arg\,min}_{v\in V_Z}J(v-u^*)$ and a comparison
function $\widetilde u=I_Zu^*\in V_Z$, linearity gives $u_Z-\widetilde u\in V_Z$.
Discrete stability and least-squares minimality yield
\begin{flalign}\label{kernel-consistency}
\|u_Z-u^*\|_{H^2(\Omega)}^2
&\lesssim \underbrace{ J(u_Z-\widetilde u)
+\|\widetilde u-u^*\|_{H^2(\Omega)}^2 
}_{\displaystyle  \text{Subtraction trick  }  u_Z-\widetilde u\in V_Z}&&\\
\makebox[12em][l]{(Consistency)}
&\lesssim J(u_Z- u^*)+J(\widetilde u-u^*)
+\|\widetilde u-u^*\|_{H^2(\Omega)}^2&&\notag\\
&\lesssim 2J(\widetilde u-u^*)
+\|\widetilde u-u^*\|_{H^2(\Omega)}^2.&&\notag
\end{flalign}
Kernel approximation theory bounds both terms on the right,
yielding convergence of least-squares kernel-based collocation. Figure~\ref{fig:three-ingredients-convergence} summarizes the roles
of regularity, stability, and consistency in this analysis.

% \medskip

\subsection{\texorpdfstring{Stable admissible neural classes}{Stable admissible neural classes}}
\begin{figure}[htbp]
\centering
\begin{tikzpicture}[
neuron/.style={
circle, draw=teal!80!black,
fill=teal!10, minimum size=5mm,
inner sep=0pt},
connection/.style={draw=gray!55, thin}]
% Use one canvas so the network and trial-class diagram share a centerline.
\begin{scope}[x=1.60cm, y=0.80cm]
\foreach \i/\y in {1/0.5,2/-0.5} \node[neuron] (I\i) at (0,\y) {};

% Hidden layers
\foreach \i/\y in {1/1,2/0,3/-1} {
\node[neuron] (H\i) at (1,\y) {};
\node[neuron] (G\i) at (2,\y) {};}

% Output layer
\node[neuron] (O) at (3,0) {};

% Connections
\foreach \i in {1,2}
\foreach \j in {1,2,3}
\draw[connection] (I\i) -- (H\j);

\foreach \i in {1,2,3}
\foreach \j in {1,2,3}
\draw[connection] (H\i) -- (G\j);

\foreach \i in {1,2,3}
\draw[connection] (G\i) -- (O);

% Labels
\node[font=\footnotesize] at (0,-1.9) {Input $x$};
\node[font=\footnotesize] at (1.5,-1.9) {Hidden layers};
\node[font=\footnotesize] at (3,-1.9) {Output $v_\theta\in\NN$};
\end{scope}

\begin{scope}[xshift=8.30cm]
\filldraw[fill=pale!65, draw=sea, line width=0.8pt]
(0,0) ellipse (2.05cm and 1.40cm);
\filldraw[fill=white, draw=warm, line width=0.9pt]
(0,-0.38) ellipse (1.60cm and 0.70cm);
\node[font=\footnotesize, align=center, text=ink] at (0,0.82)
{$\NN$\\prescribed trial class};
\node[font=\footnotesize, align=center, text=warm] at (0,-0.38)
{$\NNa$\\admissible subclass};
\end{scope}
\end{tikzpicture}
\caption{Neural trial class and its stable admissible subclass.}
\label{fig:neural-class}
\vspace{-10pt}
\end{figure}
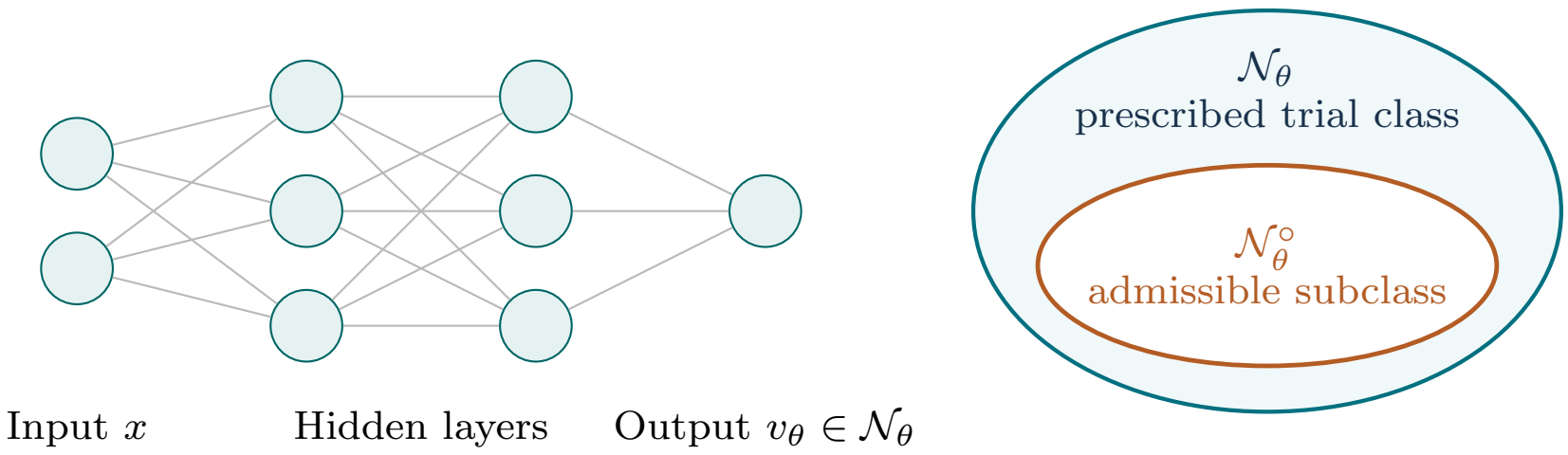

\begin{figure}[htbp]
\centering
\begin{tikzpicture}
\node[draw=warm,line width=0.5pt,rounded corners=2pt,
inner sep=5pt,minimum height=1.65cm,
text width=\dimexpr0.94\linewidth-12pt\relax,
align=left,font=\small,text=black]
{\begin{minipage}[c]{0.18\linewidth}
\centering
\textbf{\textcolor{warm}{Stability}}\\[3pt]
{\footnotesize
\hyperref[sec:bii_admissible_nn]{\textcolor{gray!75}{%
(Section~\ref*{sec:bii_admissible_nn})}}}
\end{minipage}\hspace{0.04\linewidth}%
\begin{minipage}[c]{0.8\linewidth}
\raggedright
How do we impose an appropriate BII on a trial class
that \\ is \emph{not a vector space}? Restrict to an admissible
subclass  $\NNa$\\ on which we have stability.
\end{minipage}};
\end{tikzpicture}
\par\smallskip
\begin{tikzpicture}
\node[draw=warm,line width=0.5pt,rounded corners=2pt,
  inner sep=5pt,minimum height=1.65cm,
  text width=\dimexpr0.94\linewidth-12pt\relax,
  align=left,font=\small,text=black]
{\begin{minipage}[c]{0.18\linewidth}
\centering
\textbf{\textcolor{warm}{Consistency}}\\[3pt]
{\footnotesize
\hyperref[sec:approximation-consistency]{\textcolor{gray!75}{%
(Section~\ref*{sec:approximation-consistency})}}}
\end{minipage}\hspace{0.04\linewidth}%
\begin{minipage}[c]{0.8\linewidth}
\raggedright
\setlength{\abovedisplayskip}{3pt}
\setlength{\abovedisplayshortskip}{3pt}
\setlength{\belowdisplayskip}{0pt}
\setlength{\belowdisplayshortskip}{0pt}
To use the \emph{method of comparison}, we need an accurate\\
comparison function \emph{inside the admissible subclass}:
\[
A_hu^*\in\NNa,\qquad A_hu^*\approx u^*.
\]
\end{minipage}};
\end{tikzpicture}
\caption{Stability and consistency  for minimization over the admissible subclass.}
\vspace{-10pt}
\label{fig:admissible-neural-subclass}
\end{figure}

Motivated by this kernel argument, we replace $V_Z$ by
a neural trial class. In this setting,
\begin{align}\label{eq:neural-trial-class-P}
\mathcal N_\theta
=
\{v_\theta:\theta\in\Theta_n\},
\end{align}
where $n\in\mathbb N$ indexes a prescribed family of network
architectures and $\Theta_n$ is the corresponding parameter set.
Physics-informed MLPs and KANs still minimize sampled residuals.
The trial class is different, but the stability question
remains the same, namely when small discrete loss implies
small solution error.

Fixed-architecture MLPs and KANs define trial classes that are
generally not vector spaces. Consequently, the subtraction trick
in \eqref{kernel-consistency} is no longer automatic.
In particular, $v_\theta,\widetilde v_\theta\in\NN$
does not imply $v_\theta-\widetilde v_\theta\in\NN$.
For comparison, consider the kernel solution $u_Z\in V_Z$.
Applying the regularity estimate \eqref{eq:bdyreg}
and the sampling inequalities~\cite{CheungLingSchaback2018H2}
to the actual error $u_Z-u^*$ gives
{\small
\medmuskip=2mu
\thickmuskip=3mu
\begin{align}
\|u_Z-u^*\|_{H^2(\Omega)}^2
\lesssim
J(u_Z-u^*)
+\underbrace{
h_X^{2m-4}\|\mathcal L(u_Z-u^*)\|_{H^{m-2}(\Omega)}^2
+h_Y^{2m-4}\|u_Z-u^*\|_{H^{m-\frac12}(\partial\Omega)}^2
}_{\text{interior \& boundary strong-norm remainder}}.\notag
\end{align}}%
\noindent Operator continuity and the trace theorem then bound
both strong norms by $\norm{u_Z-u^*}_{H^m(\Omega)}$.
Since $u^*$ need not belong to $V_Z$,
\begin{align*}
\norm{u_Z-u^*}_{H^m}
\le
\underbrace{\norm{u_Z-I_Zu^*}_{H^m}}_{\text{BII in }V_Z}
+
\underbrace{\norm{I_Zu^*-u^*}_{H^m}}_{\text{Interp. theory}}.
\end{align*}
We then apply BII within the trial space and use oversampling
to absorb the error-dependent part of the remainder.

For the neural route, general neural network classes $\NN$
do not automatically satisfy the BII available on kernel
trial spaces. This leads to the central problem of designing
a stable admissible neural class $\NNa \subset \NN$,  as illustrated in
Figure~\ref{fig:neural-class}. 
Rather than applying BII to the actual error $v_\theta-u^*$,
we impose a PDE-dependent residual-BII constraint directly
on the interior and boundary residuals to define an admissible
subclass $\NNa$.
For any trial network $v_\theta\in\NN$, define
\begin{align}
\begin{aligned}
r_\Omega(v_\theta):=\mathcal L(v_\theta-u^*)=\mathcal Lv_\theta-f \quad \text{and} \quad 
r_{\partial\Omega}(v_\theta):=(v_\theta-u^*)|_{\partial\Omega}
=v_\theta|_{\partial\Omega}-g.
\end{aligned}
\end{align}
The regularity--sampling estimate becomes
% {\small
% \medmuskip=2mu
% \thickmuskip=3mu
\begin{align}
  % \label{eq:standard-sampling-remainder}
\|v_\theta- u^*\|_{H^2(\Omega)}^2
\lesssim 
J(v_\theta-u^*) 
+h_X^{2m-4}\|r_\Omega(v_\theta)\|_{H^{m-2}(\Omega)}^2
+h_Y^{2m-4}
\|r_{\partial\Omega}(v_\theta)\|_{H^{m-\frac{1}{2}}(\partial\Omega)}^2. \notag
\end{align}
By operator continuity and the trace estimate, the lower-order
residual norms satisfy
% \begin{align}
% \alert{\norm{r_\Omega(v_\theta)}_{L^2(\Omega)}^2}\lesssim  {\norm{v_\theta-u^*}_{H^2(\Omega)}^2}, \quad 
% \alert{\norm{r_{\partial\Omega}(v_\theta)}_{H^{3/2}(\partial\Omega)}^2} \lesssim  {\norm{v_\theta-u^*}_{H^2(\Omega)}^2}. 
% \end{align}
\begin{align}
\|r_\Omega(v_\theta)\|_{L^2(\Omega)}^2
\le C_{\mathcal L}\|v_\theta-u^*\|_{H^2(\Omega)}^2, \,\,
\|r_{\partial\Omega}(v_\theta)\|_{H^{3/2}(\partial\Omega)}^2
\le C_{\mathrm{tr}}\|v_\theta-u^*\|_{H^2(\Omega)}^2.
\label{continuity-L}
\end{align}
Thus, to control the strong-norm remainders from sampling,
we only need the following residual-BII bounds, scaled by
the design scale $h_\theta$, to hold on the admissible
subclass $\NNa$:
\begin{align}\label{bii-domain}
\|r_\Omega(v_\theta)\|_{H^{m-2}(\Omega)}^2
&\le C_\Omega h_\theta^{-(2m-4)}
\|r_\Omega(v_\theta)\|_{L^2(\Omega)}^2,\\
\|r_{\partial\Omega}(v_\theta)\|_{H^{m-\frac12}(\partial\Omega)}^2
&\le C_{\partial\Omega}h_\theta^{-(2m-4)}
\|r_{\partial\Omega}(v_\theta)\|_{H^{3/2}(\partial\Omega)}^2.
\label{bii-boundary}
\end{align}
The computational bottleneck is that demanding these full
inverse estimates directly is not computationally friendly. They involve fractional-order
Sobolev norms and large derivative gaps, such as the gap
of $m-2$ orders between $H^{m-2}(\Omega)$ and $L^2(\Omega)$.

We therefore provide a tractable fix through a practical
prescription based on a single-gap BII.
For $D\in\{\Omega,\partial\Omega\}$, choose
$C_{\mathrm{BII},D}>0$ and $\beta_D\in\mathbb N$,
and prescribe a small fixed gap $\varepsilon_*>0$:
\begin{equation}\label{eq:intro-residual-bii}
\|r_D(v_\theta)\|_{H^{\beta_D}(D)}
\le C_{\mathrm{BII},D}h_\theta^{-\varepsilon_*}
\|r_D(v_\theta)\|_{H^{\beta_D-\varepsilon_*}(D)},
\end{equation}
where the BII constants are independent of the network
design scale $0<h_\theta\le1$.
Choosing an integer gap, such as $\varepsilon_*=1$,
also avoids fractional-order norms in these constraints.
We define the BII-admissible neural subclass $\NNa$
through the single-gap residual-BII condition
\eqref{eq:intro-residual-bii}.

The residual-BII prescription above specifies the admissible
subclass $\NNa$ for the stability task in
Figure~\ref{fig:admissible-neural-subclass}.
For sufficiently dense residual sampling relative to $h_\theta$,
the conditions of Section~\ref{sec:bii_admissible_nn} yield
\begin{flalign}\label{eq:intro-neural-stability}
\makebox[12em][l]{(Stability)} & \|v_\theta-u^*\|_{H^2(\Omega)}^2
\lesssim J(v_\theta-u^*),
\qquad v_\theta\in\NNa,&&
\end{flalign}
where the constant is independent of $h_\theta,h_X,h_Y$.
The canonical objective $J$ is further extended to a family
of weighted objectives $J_\gamma$, with corresponding
discrete stability estimates.
Using this weighted family $J_\gamma$, Section~3 introduces
physics-balanced unconstrained (PBU) and constrained (PBC)
minimization on general trial classes.
PBU minimizes $J_\gamma$ over the full trial class,
whereas PBC minimizes the same objective over its
BII-admissible subclass.

Section~\ref{sec:approximation-consistency} addresses the
consistency and convergence in Figure~\ref{fig:admissible-neural-subclass}.
It establishes convergence of PBC and PBU
on general trial classes, together with $H^2$-convergence rates.
For the prescribed $\tanh$ network family, we fix arbitrary
positive BII thresholds and construct approximants
$A_hu^*\to u^*$  as $h_\theta\to0$ such that
\begin{flalign}\label{eq:intro-neural-consistency}
\makebox[12em][l]{(Consistency)}&
A_hu^*\in\NNa,
\qquad J_\gamma(A_hu^*-u^*)\to0. &&
\end{flalign}
This provides a constructive consistency proof within
the prescribed admissible neural subclass.
For computation, the continuous BII constraints are replaced
by empirical soft BII penalties.
The numerical experiments in
Section~\ref{sec:numerical-experiments} compare the two
formulations and examine how sampling refinement,
network refinement, and the BII thresholds affect the error.

\medskip

\section{Admissible neural classes and discrete stability}
\label{sec:bii_admissible_nn}
Uniform residual BII bounds define an admissible subclass of
the prescribed neural network class. We establish discrete stability on this
subclass under the stated assumptions, with sampling requirements depending
on the design scale $h_\theta$ and the BII constants.

\subsection{A single-gap BII and the admissible neural subclass}
\label{subsec:bii-admissible-neural-residuals}
% To control the strong-norm remainders from sampling, we only need BII on the residuals.
To control the strong-norm remainders in the sampling estimates, we impose single-gap BII on the interior and boundary residuals.
% For statements common to both residuals, 
% $r_D(v_\theta)$ denote the corresponding domain and residual.

\begin{definition}[BII-admissible neural subclass $\mathcal N_\theta^\circ$]
\label{ass:endpoint-bii}
Choose $C_{\mathrm{BII},D}>0$ and $\beta_D\in\mathbb N$  for $D\in\{\Omega,\partial\Omega\}$, and prescribe a small fixed gap $0<\varepsilon_*\le1$. Consider the following single-gap BII condition:
\begin{align}\label{BII}
\|r_D(v_\theta)\|_{H^{\beta_D}(D)}
\le C_{\mathrm{BII},D}h_\theta^{-\varepsilon_*}
\|r_D(v_\theta)\|_{H^{\beta_D-\varepsilon_*}(D)},
\quad D\in\{\Omega,\partial\Omega\},
\end{align}
where $C_{\mathrm{BII},D}$ is independent of network design scale $h_\theta$ with $0<h_\theta\le1$.
We define the admissible subclass $\mathcal N_\theta^\circ\subseteq\mathcal N_\theta$ purely through this condition:
\begin{equation}\label{eq:bii_neural_family}
\mathcal N_\theta^\circ
:=
\left\{
v_\theta\in\mathcal N_\theta:
\begin{array}{l}
r_D(v_\theta)\in H^{\beta_D}(D)
\text{ and satisfies the BII \eqref{BII}}\\
\text{for }D\in\{\Omega,\partial\Omega\}
\end{array}
\right\}.
\end{equation}
% Here the PDE, $\varepsilon_*$, and the BII constants are fixed as stated above.
\end{definition}

For fixed
$\mathcal N_\theta$ and $h_\theta$, changing $f$ or $g$ changes the residuals
and may therefore change the admissible class.  Changing the coefficients of
$\mathcal L$ has the same effect.  Thus the admissible class depends on the PDE
data.  
% We use $\mathcal N_\theta^\circ$ as an abbreviation for the full notation $\mathcal N_\theta^\circ(f,g;h_\theta,\beta_\Omega,\beta_{\partial\Omega})$. 
% Throughout this section, we assume Sobolev interpolation
% at the prescribed BII orders.
% Figure~\ref{fig:admissible-neural-subclass} illustrates the inclusion
% $\NNa\subseteq\NN$.
The following lemma derives full inverse scale estimates from the $\varepsilon_*$-gap residual BII.

\begin{lemma}[Residual inverse inequality]
\label{thm:abstract-bernstein-extension}
In the setting of Definition~\ref{ass:endpoint-bii}, for any
$v_\theta\in\mathcal N_\theta^\circ$ and $D\in\{\Omega,\partial\Omega\}$,
the single gap $\varepsilon_*$ at order $\beta_D$ automatically recovers
the full inverse bound down to the base norms:
% \begin{equation}\label{eq:abstract-bernstein-full}
% \begin{gathered}
% \|r_D(v_\theta)\|_{H^\alpha(D)}
% \le C_{D,\beta_D,\varepsilon_*,\alpha,t}h_\theta^{-(\alpha-t)}
% \|r_D(v_\theta)\|_{H^t(D)},\\
% 0\le t\le\beta_D-\varepsilon_*,\qquad t\le\alpha\le\beta_D.
% \end{gathered}
% \end{equation}
\begin{equation}\label{eq:abstract-bernstein-full}
\begin{gathered}
\|r_D(v_\theta)\|_{H^\alpha(D)}
\le C h_\theta^{-(\alpha-t)}
\|r_D(v_\theta)\|_{H^t(D)}, \quad
0\le t\le\beta_D-\varepsilon_*, \,  t\le\alpha\le\beta_D,
\end{gathered}
\end{equation}
where $C$ depends on $D$, the BII parameters and interpolation constants,
and is independent of the network design scale $h_\theta$.
% and the network architecture.
% where the constant $C$ may depend on $D$, $\beta_D$, $\varepsilon_*$, $\alpha$, $t$,
% $C_{\mathrm{BII},D}$, and the Sobolev interpolation constants.  It is
% uniform along architecture refinement and independent of discretization scale
% $h_\theta$.
\end{lemma}
% \medskip

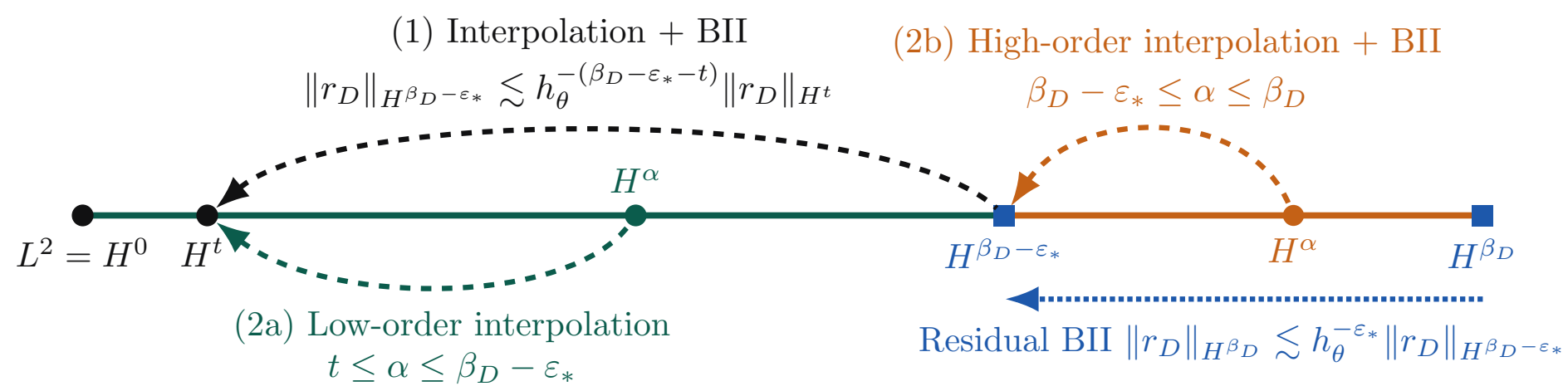
\begin{figure}[t]
\centering
\definecolor{BIIblue}{HTML}{1d58ab}
\definecolor{LowGreen}{HTML}{0b5c4c}
\definecolor{HighOrange}{HTML}{c46116}
\definecolor{InkGray}{HTML}{111112}
\resizebox{.94\linewidth}{!}{%
\begin{tikzpicture}[
>=Latex,
% every node/.style={font=\small},
% title/.style={font=\small},
% tinylabel/.style={font=\small,align=center},
tinylabel/.style={align=center},
flowarc/.style={<-,dashed,line width=1.5pt,
  shorten <=3pt,shorten >=3pt}]
\path[use as bounding box] (-7.65,0.25) rectangle (7.65,4.18);
\begin{scope}[yshift=-1.5mm]

\coordinate (L2)    at (-7.0,2.10);
\coordinate (Ht)    at (-5.75,2.10);
\coordinate (Ht')    at (-5.8,2.10);
\coordinate (Hlow)  at (-1.45,2.10);
\coordinate (Hpiv)  at ( 2.25,2.10);
\coordinate (Hhigh) at ( 5.15,2.10);
\coordinate (Htop)  at ( 7.05,2.10);

\draw[InkGray,line width=1pt] (L2)--(Htop);
\draw[LowGreen,line width=1.8pt] (L2)--(Hpiv);
\draw[HighOrange,line width=1.8pt] (Hpiv)--(Htop);

\fill[InkGray] (L2) circle (3.1pt);
\fill[InkGray] (Ht) circle (3.1pt);
\fill[LowGreen] (Hlow) circle (3.1pt);
\fill[BIIblue] ($(Hpiv)+(-3pt,-3pt)$) rectangle
($(Hpiv)+( 3pt, 3pt)$);
\fill[HighOrange] (Hhigh) circle (3pt);
\fill[BIIblue] ($(Htop)+(-3pt,-3pt)$) rectangle
($(Htop)+(3pt, 3pt)$);

\node[below=3pt,InkGray] at (L2) {$L^2=H^0$};
\node[below=3pt,InkGray] at (Ht') {$H^t$};
\node[above=3pt,LowGreen] at (Hlow) {$H^\alpha$};
\node[below=3pt,BIIblue] at (Hpiv) {$H^{\beta_D-\varepsilon_*}$};
\node[below=3pt,HighOrange] at (Hhigh) {$H^\alpha$};
\node[below=3pt,BIIblue] at (Htop) {$H^{\beta_D}$};

\draw[flowarc,InkGray]
(Ht) .. controls (-4.45,3.20) and (1.05,3.20) .. (Hpiv)
node[pos=.46,above=3pt,tinylabel,InkGray]
{(1) Interpolation $+$ BII\\[2pt]
$\|r_D\|_{H^{\beta_D-\varepsilon_*}}
\lesssim h_\theta^{-(\beta_D-\varepsilon_*-t)}\|r_D\|_{H^t}$};

\draw[flowarc,LowGreen]
(Ht) .. controls (-4.55,1.18) and (-2.10,1.18) .. (Hlow)
node[pos=.52,below=3pt,tinylabel,LowGreen]
{(2a) Low-order interpolation\\
$t\le\alpha\le \beta_D-\varepsilon_*$};

\draw[flowarc,HighOrange]
(Hpiv) .. controls (3.00,3.20) and (4.70,3.20) .. (Hhigh)
node[pos=.52,above=3.1pt,tinylabel,HighOrange]
{(2b) High-order interpolation $+$ BII\\[2pt]
$\beta_D-\varepsilon_*\le\alpha\le \beta_D$};

\draw[->,BIIblue,densely dotted,line width=1.5pt]
($(Htop)+(0,-.84)$)--($(Hpiv)+(0,-.84)$)
node[midway,below=2pt,tinylabel,BIIblue]
{Residual BII
$\|r_D\|_{H^{\beta_D}}
\lesssim h_\theta^{-\varepsilon_*}
\|r_D\|_{H^{\beta_D-\varepsilon_*}}$};
\end{scope}
\end{tikzpicture}
}
\caption{$\varepsilon_*$-gap generates full inverse scale.
Interpolation and the residual BII first yield the estimate at order
$\beta_D-\varepsilon_*$ in step (1). This estimate then extends to lower
orders by interpolation in (2a), and to higher orders by interpolation
and one further application of the BII in (2b).}
\label{fig:single-gap-bii-proof}
\vspace{-15pt}
\end{figure}

\begin{proof}
Fix $D$ and $0\le t\le \beta_D-\varepsilon_*$, and set
$\vartheta=(\beta_D-\varepsilon_*-t)/(\beta_D-t)$.
As in step (1) of Figure~\ref{fig:single-gap-bii-proof},
interpolating between $H^t(D)$ and $H^{\beta_D}(D)$,
applying the BII condition \eqref{BII}
to the higher-order norm, and rearranging the resulting
inequality (the zero case being immediate) gives
\begin{equation}\label{bii-step1}
\|r_D(v_\theta)\|_{H^{\beta_D-\varepsilon_*}(D)}
\le C h_\theta^{-(\beta_D-\varepsilon_*-t)}
\|r_D(v_\theta)\|_{H^t(D)}.
\end{equation}
Here we used
$\varepsilon_*\vartheta/(1-\vartheta)
=\beta_D-\varepsilon_*-t$.

In step (2a), we consider
$t\le\alpha\le\beta_D-\varepsilon_*$.
Interpolation between $H^t(D)$ and
$H^{\beta_D-\varepsilon_*}(D)$, together with
\eqref{bii-step1}, gives
\eqref{eq:abstract-bernstein-full}. 
In step (2b), we consider
$\beta_D-\varepsilon_*\le\alpha\le\beta_D$.
We interpolate between the
$H^{\beta_D-\varepsilon_*}(D)$ and $H^{\beta_D}(D)$
norms and apply BII once more.
This supplies the factor
$h_\theta^{-(\alpha-\beta_D+\varepsilon_*)}$.
Combining this with \eqref{bii-step1} gives
$h_\theta^{-(\alpha-t)}$ and hence
\eqref{eq:abstract-bernstein-full}.
The endpoint cases follow directly.
\end{proof}

\subsection{Fully discrete $H^2$ stability on $\NNa$}
\label{sec:WLS-convergence}
We first establish discrete stability on $\mathcal N_\theta^\circ$
under suitable oversampling and then show that a lower-order boundary BII
suffices for the same estimate.
\begin{samepage}
\begin{lemma}[Discrete stability in canonical form for all $v_\theta\in\NNa$]
\label{thm:standard-residual-stability}
Let $u^*\in H^m(\Omega)$ solve \eqref{eq:PDE} with $m>2+d/2$.
For any admissible class $\mathcal N_\theta^\circ\subset H^m(\Omega)$ defined in
Definition~\ref{ass:endpoint-bii} with BII orders
$\beta_\Omega\ge m-2$ and 
$\beta_{\partial\Omega}\ge m-\frac12$,
let $0<h_X\le h_Y<1$ be sufficiently small and satisfy
the oversampling condition
\begin{equation}\label{eq:standard-density-conditions}
C_{\rm abs}
\left[
\left(\frac{h_X}{h_\theta}\right)^{2m-4}
+
\left(\frac{h_Y}{h_\theta}\right)^{2m-4}
\right]\le\frac12.
\end{equation}
Then every $v_\theta\in\mathcal N_\theta^\circ$ satisfies
\begin{equation}\label{eq:canonical-residual-stability}
\|v_\theta-u^*\|_{H^2(\Omega)}^2
\lesssim  J(v_\theta-u^*):= \,
h_X^d\|r_\Omega(v_\theta)\|_{\ell^2(X)}^2
+
h_Y^{d-4}
\|r_{\partial\Omega}(v_\theta)\|_{\ell^2(Y)}^2.
\end{equation}
where the constants 
are independent
of $h_X,h_Y$ and the network design scale $h_\theta$.
\end{lemma}
\end{samepage}

\begin{proof}
Recall the regularity--sampling estimate. 
{\small
\begin{align}
\label{eq:standard-sampling-remainder}
\|v_\theta- u^*\|_{H^2(\Omega)}^2
\lesssim 
J(v_\theta-u^*) 
+h_X^{2m-4}\|r_\Omega(v_\theta)\|_{H^{m-2}(\Omega)}^2
+h_Y^{2m-4}
\|r_{\partial\Omega}(v_\theta)\|_{H^{m-\frac{1}{2}}(\partial\Omega)}^2.
\end{align}}%
We control its strong-norm remainders using the residual
BII estimates in Lemma~\ref{thm:abstract-bernstein-extension}.
Under the prescribed BII orders, Lemma~\ref{thm:abstract-bernstein-extension} yields the full inverse bounds
on $\NNa$
\begin{align}
\|r_\Omega(v_\theta)\|_{H^{m-2}(\Omega)}^2
&\le C_\Omega h_\theta^{-(2m-4)}
\|r_\Omega(v_\theta)\|_{L^2(\Omega)}^2,\\
\|r_{\partial\Omega}(v_\theta)\|_{H^{m-\frac12}(\partial\Omega)}^2
&\le C_{\partial\Omega}h_\theta^{-(2m-4)}
\|r_{\partial\Omega}(v_\theta)\|_{H^{3/2}(\partial\Omega)}^2.
\end{align}
% Operator continuity and the trace theorem give
% \begin{align}
% \|r_\Omega(v_\theta)\|_{L^2(\Omega)}^2
% \le C_{\mathcal L}\|v_\theta-u^*\|_{H^2(\Omega)}^2, \,\,
% \|r_{\partial\Omega}(v_\theta)\|_{H^{3/2}(\partial\Omega)}^2
% \le C_{\mathrm{tr}}\|v_\theta-u^*\|_{H^2(\Omega)}^2.
% \label{continuity-L}
% \end{align}
For the interior remainder, combining \eqref{bii-domain}
with the continuity estimate \eqref{continuity-L} gives
\begin{align}\label{eq:standard-interior-remainder}
h_X^{2m-4}
\|r_\Omega(v_\theta)\|_{H^{m-2}(\Omega)}^2
&\le
C\left(\frac{h_X}{h_\theta}\right)^{2m-4}
\|r_\Omega(v_\theta)\|_{L^2(\Omega)}^2 \\
&\le
C\left(\frac{h_X}{h_\theta}\right)^{2m-4}
\|v_\theta-u^*\|_{H^2(\Omega)}^2.
\notag
\end{align}
For the boundary remainder, combining \eqref{bii-boundary}
with the trace estimate \eqref{continuity-L} gives
\begin{align}\label{eq:standard-boundary-remainder}
h_Y^{2m-4}
\|r_{\partial\Omega}(v_\theta)\|_
{H^{m-\frac12}(\partial\Omega)}^2
&\le
C\left(\frac{h_Y}{h_\theta}\right)^{2m-4}
\|r_{\partial\Omega}(v_\theta)\|_
{H^{3/2}(\partial\Omega)}^2 \\
&\le
C\left(\frac{h_Y}{h_\theta}\right)^{2m-4}
\|v_\theta-u^*\|_{H^2(\Omega)}^2.
\notag
\end{align}
Substituting \eqref{eq:standard-interior-remainder} and
\eqref{eq:standard-boundary-remainder} into
\eqref{eq:standard-sampling-remainder} and collecting
constants yields
\begin{align*}
\|v_\theta-u^*\|_{H^2(\Omega)}^2
\le{}& CJ(v_\theta-u^*)+C_{\rm abs}
\left[
\left(\frac{h_X}{h_\theta}\right)^{2m-4}
+\left(\frac{h_Y}{h_\theta}\right)^{2m-4}
\right]
\|v_\theta-u^*\|_{H^2(\Omega)}^2.
\end{align*}
By the oversampling condition
\eqref{eq:standard-density-conditions}, the coefficient
of $\|v_\theta-u^*\|_{H^2(\Omega)}^2$ in the last term
is at most $1/2$.
Absorbing this term into the left-hand side yields
\eqref{eq:canonical-residual-stability}.
Throughout the paper, for $d=1$, take $Y=\partial\Omega$, interpret all boundary Sobolev
norms as the endpoint $\ell^2$-norm, and replace the boundary sampling
assumptions by exact endpoint evaluation, with $h_Y=h_X$ serving as an
auxiliary positive scale. 
\end{proof}

\begin{table}[t]
\centering
\caption{Parameter prescription for the residual-BII with small
$\varepsilon_*>0$}
\label{tab:residual-bii-parameters}
\renewcommand{\arraystretch}{1.18}
\setlength{\tabcolsep}{9pt}
\begin{tabular}{@{}lccl@{}}
\toprule
Domain $D$ & Sampling order $\alpha$
& Stability target $t$ & Data regularity \\
\midrule
Interior $\Omega$ & $>d/2$ & $0$
& $f\in H^\alpha(\Omega)$ \\
\addlinespace[3pt]
Boundary $\partial\Omega$ & $>\frac12\max\{d,3\}$
& $\tfrac32$ & $g\in H^\alpha(\partial\Omega)$ \\
\bottomrule
\end{tabular}
\end{table}

Table~\ref{tab:residual-bii-parameters} summarizes the sampling orders,
stability targets, and data regularity used in the preceding
sampling--BII argument.
For the direct application of
Lemma~\ref{thm:abstract-bernstein-extension} at the stability target $t$,
the BII orders must satisfy
$\alpha\le\beta_D$ and $0\le t\le\beta_D-\varepsilon_*$.

To simplify computation in practice, we choose integer BII orders
$\beta_D$ with integer gap $\varepsilon_*=1$, avoiding fractional
norms in the constraints.
We further seek to reduce the derivative orders required for stability.
With the boundary target $t=3/2$, this direct argument requires at least
an $H^3$ to $H^2$ BII.
The following lemma shows that the $H^2$ to $H^1$ boundary BII already
suffices for the same canonical stability estimate.

\begin{lemma}[Discrete stability under lower-order BII
for all $v_\theta\in\NNa$]
\label{thm:int-residual}
Let $u^*\in H^m(\Omega)$ solve \eqref{eq:PDE}
with $m>2+d/2$ and $d\le3$.
For any BII-admissible class
$\mathcal N_\theta^\circ\subset H^2(\Omega)$ with orders
$\beta_\Omega=\beta_{\partial\Omega}=2$ and gap $\varepsilon_*=1$,
let $0<h_X\le h_Y<1$ be sufficiently small and satisfy
the oversampling condition
\begin{equation}\label{eq:density-conditions}
\mathcal C_{\rm abs}'
\left[
\left(\frac{h_X}{h_\theta}\right)^4
+
\left(\frac{h_Y}{h_\theta}\right)^3
\right]\le\frac12.
\end{equation}
Then every $v_\theta\in\mathcal N_\theta^\circ$ satisfies
\begin{equation*}
\|v_\theta-u^*\|_{H^2(\Omega)}^2
\lesssim J(v_\theta-u^*).
\end{equation*}
% and, consequently, the weighted estimate
% \eqref{eq:parameterized-discrete-stability}
% for every $\gamma\in[0,2]$.
% The constants are independent of $h_X,h_Y,h_\theta$
% and the particular network.
\end{lemma}

% \begin{theorem}[Discrete stability under a lower-order boundary BII]
% \label{thm:int-residual}
% Let $d=2,3$, and let $u^*\in H^2(\Omega)$ be the solution of
% \eqref{eq:PDE}. Assume the elliptic $H^2$ regularity estimate,
% operator and trace bounds used in Section~\ref{sec:WLS-convergence},
% and the second-order residual sampling and Sobolev interpolation
% estimates. Use Definition~\ref{ass:endpoint-bii} with
% $\mathcal N_\theta\subset H^2(\Omega)$ and
% $\beta_\Omega=\beta_{\partial\Omega}=2$, $\varepsilon_*=1$.
% Let $0<h_X\le h_Y<1$ be sufficiently small and satisfy
% \begin{equation}\label{eq:density-conditions}
% \mathcal C_{\rm abs}^{\prime}
% \left[
% \left(\frac{h_X}{h_\theta}\right)^4
% +\left(\frac{h_Y}{h_\theta}\right)^4
% \right]\le\frac12,
% \end{equation}
% where $\mathcal C_{\rm abs}^{\prime}\ge1$ is a sufficiently large fixed constant.
% Then every $v_\theta\in\mathcal N_\theta^\circ$ satisfies the canonical
% estimate \eqref{eq:canonical-residual-stability} and, for every
% $\gamma\in[0,2]$, the weighted estimate
% \eqref{eq:parameterized-discrete-stability}, with the same definitions
% of $W_\gamma$ and $J_\gamma$.
% where $C$ is independent of $h_X,h_Y$ and the network design scale $h_\theta$.
% \end{theorem}

\begin{proof}
Recall from the regularity estimate
\eqref{eq:bdyreg} with $k=0$ that
\begin{equation*}
\|v_\theta-u^*\|_{H^2(\Omega)}^2
\lesssim
\|r_\Omega(v_\theta)\|_{L^2(\Omega)}^2
+
\|r_{\partial\Omega}(v_\theta)\|_{H^{3/2}(\partial\Omega)}^2, \quad v_\theta\in\mathcal N_\theta^\circ.
\end{equation*}
For the interior residual, we follow the sampling--BII argument
in the preceding proof with $\beta_\Omega=2>d/2$ to obtain
\begin{align}\label{eq:low-order-interior-control}
\|r_\Omega(v_\theta)\|_{L^2(\Omega)}^2
&\le C h_X^d\|r_\Omega(v_\theta)\|_{\ell^2(X)}^2+C\left(\frac{h_X}{h_\theta}\right)^4
\|v_\theta-u^*\|_{H^2(\Omega)}^2.
\end{align}

For the boundary residual, we reverse the order.
The residual inverse estimate in
Lemma~\ref{thm:abstract-bernstein-extension} with
$t=0$ and $\alpha=3/2$, followed by the boundary sampling estimate
\cite[Lemma~3.1]{CheungLingSchaback2018H2}
for $v_\theta-u^*\in H^2(\Omega)$, gives
\begin{align}\label{eq:low-order-boundary-control}
\|r_{\partial\Omega}(v_\theta)\|_{H^{3/2}(\partial\Omega)}^2
&\le C h_\theta^{-3}
\|r_{\partial\Omega}(v_\theta)\|_{L^2(\partial\Omega)}^2\notag\\
&\le C h_\theta^{-3}\Bigl[
h_Y^3\|v_\theta-u^*\|_{H^2(\Omega)}^2
+h_Y^{d-1}\|r_{\partial\Omega}(v_\theta)\|_{\ell^2(Y)}^2
\Bigr]\notag\\
&\le C\left(\frac{h_Y}{h_\theta}\right)^3
\|v_\theta-u^*\|_{H^2(\Omega)}^2
+C h_Y^{d-4}\|r_{\partial\Omega}(v_\theta)\|_{\ell^2(Y)}^2.
\end{align}
Combining these estimates yields, 
\begin{align}
\|v_\theta-u^*\|_{H^2(\Omega)}^2
\le C J(v_\theta-u^*)
+\mathcal C_{\rm abs}'
\left[
\left(\frac{h_X}{h_\theta}\right)^4
+\left(\frac{h_Y}{h_\theta}\right)^3
\right]\|v_\theta-u^*\|_{H^2(\Omega)}^2.
\end{align}
Absorbing the last term using \eqref{eq:density-conditions}
proves \eqref{eq:canonical-residual-stability}.
\end{proof}

% Here $M_X=|X|$, $M_Y=|Y|$, and $M=M_X+M_Y$.
% \par canonical estimate \eqref{eq:canonical-residual-stability}

\subsection{Generalized weighted least squares}
The preceding subsection establishes canonical discrete stability \eqref{eq:canonical-residual-stability} on
$\NNa$. We now derive a family of weighted stability estimates.

% Starting from the canonical form, we obtain the following discrete stability family. 
\begin{theorem}[Discrete stability family for all $v_\theta\in\NNa$]
\label{cor:physics-balanced-stability}
Under the assumptions of either
Lemma~\ref{thm:standard-residual-stability}
or Lemma~\ref{thm:int-residual}, the following estimate holds for all
$\gamma\in[0,2]$ and $v_\theta\in\NNa$
\begin{align}\label{eq:parameterized-discrete-stability}
\|v_\theta-u^*\|_{H^2(\Omega)}^2
&\lesssim J_\gamma(v_\theta-u^*)\\
&:= h_X^{\frac d2\gamma}
h_Y^{d-4+2\gamma-\frac d2\gamma}
\|r_\Omega(v_\theta)\|_{\ell^2(X)}^2
+h_Y^{d-4}
\|r_{\partial\Omega}(v_\theta)\|_{\ell^2(Y)}^2. \notag
\end{align}
\end{theorem}

\begin{proof}
By \cite[Lemma~4.1]{CheungLingSchaback2018H2},
for $a,b>0$, $0<\epsilon<1$, and $0\le\gamma\le2$,
we have
\begin{equation}\label{family-gamma}
(\epsilon a+b)^2
\le 2\left(\epsilon^\gamma a^2+b^2\right).
\end{equation}
Set $\varepsilon_{XY}=(h_X/h_Y)^{d/2}h_Y^2$.
Since $0<h_X\le h_Y<1$, we have
$0<\varepsilon_{XY}<1$.
Combining the canonical estimate
\eqref{eq:canonical-residual-stability} with \eqref{family-gamma},
applied with
$\epsilon=\varepsilon_{XY}$,
$a=\|r_\Omega(v_\theta)\|_{\ell^2(X)}$, and
$b=\|r_{\partial\Omega}(v_\theta)\|_{\ell^2(Y)}$,
we obtain, for $0\le\gamma\le2$,
\begin{align}
\|v_\theta-u^*\|_{H^2(\Omega)}^2
&\le C h_Y^{d-4}\Bigl(
\varepsilon_{XY}\|r_\Omega(v_\theta)\|_{\ell^2(X)}
+\|r_{\partial\Omega}(v_\theta)\|_{\ell^2(Y)}
\Bigr)^2 \notag\\
&\le 2C\Bigl(
h_X^{\frac d2\gamma}
h_Y^{d-4+2\gamma-\frac d2\gamma}
\|r_\Omega(v_\theta)\|_{\ell^2(X)}^2
+h_Y^{d-4}
\|r_{\partial\Omega}(v_\theta)\|_{\ell^2(Y)}^2
\Bigr).
\label{eq:expanded-weighted-stability}
\end{align}
This proves
\eqref{eq:parameterized-discrete-stability}.
\end{proof}

\medskip
% The residual expression in
% \eqref{eq:parameterized-discrete-stability} defines
% $J_\gamma(v_\theta-u^*)$ for every
% $v_\theta\in\mathcal N_\theta$.
% The stability estimate holds on the admissible subclass
% $\NNa$.
% To express this functional in weighted least-squares form,

\section{Physics-balanced formulations on general trial classes}
% The residual-BII arguments and discrete stability estimates of
% Section~\ref{sec:bii_admissible_nn} apply to general trial classes under
% the same regularity, BII, and sampling assumptions.

Let $\{\NN\}_{h_\theta}$ be a family of general
trial classes with assigned design scales $h_\theta$.
For fixed residual-BII orders $\beta_D$ high enough for stability
and a fixed gap $\varepsilon_*$, each threshold vector
$\mathbf C_{\mathrm{BII}}=(C_{\mathrm{BII},\Omega},C_{\mathrm{BII},\partial\Omega})$
independent of $h_\theta$ determines a corresponding family
$\{\NNa\}_{h_\theta}$ of admissible subclasses, with
$\NNa\subseteq\NN$ at every scale.
The same analysis also applies to essentially bounded sequences
$C_{\mathrm{BII},D}(h_\theta)$.
The size of these thresholds determines the size of $\NNa$ and governs
the oversampling requirement needed for stability. 
% A stability restriction is useless if it eliminates all accurate networks.

Under the corresponding sampling conditions, $\NNa$ inherits the discrete
stability estimates established in Section~\ref{sec:bii_admissible_nn}. In this section, we formulate PBC minimization on $\NNa$ and PBU
minimization on the full trial class $\NN$.
\subsection{PBC and PBU minimization problems} 

% This subsection will use the generalized weighted least squares $J_\gamma$ to formulate the PBC and PBU minimization problems.
Since the loss directly bounds the error $\|v_\theta - u^*\|_{H^2(\Omega)}^2 \lesssim J_\gamma(v_\theta-u^*)$ in $\NNa$, we define the physics-balanced
constrained (PBC) and physics-balanced unconstrained (PBU) solutions by
\begin{align}
u_\theta^{\mathrm{PBC}_\gamma}
&\in\argmin_{v_\theta\in\NNa}
J_\gamma(v_\theta-u^*) \qquad \text{and}\qquad
\label{eq:pbc-minimization-problem}
u_\theta^{\mathrm{PBU}_\gamma}
\in\argmin_{v_\theta\in\NN}
J_\gamma(v_\theta-u^*),
% \label{eq:pb-minimization-problem}
\end{align}
whenever these minimizers exist. Here, PBC refers to minimization over a stable
subclass (e.g., with a BII constraint), whereas PBU
refers to minimization over the full, unconstrained trial class.

Define the physics-balanced boundary weight
$W_\gamma$ for $\gamma\in[0,2]$ by
\begin{align}\label{eq:boundary-weight-family}
W_\gamma:=
\left(\frac{h_Y}{h_X}\right)^{\frac d2\gamma}h_Y^{-2\gamma}.
\end{align}
Then $J_\gamma$ can be written as
\begin{align}
J_\gamma(v_\theta-u^*)
=h_Y^{d-4}W_\gamma^{-1}
\left(
\|r_\Omega(v_\theta)\|_{\ell^2(X)}^2
+W_\gamma
\|r_{\partial\Omega}(v_\theta)\|_{\ell^2(Y)}^2
\right).
\end{align}
The endpoint choices of $\gamma$ recover two familiar
residual objectives:
\begin{enumerate}[
label=(\roman*),
leftmargin=2em,
labelsep=0.4em,
itemsep=3pt,
topsep=3pt
]
\item $\gamma=2$ (optimal WLS) recovers the scale-balanced
functional in \eqref{eq:canonical-residual-stability}:
\begin{align}
J_2(v_\theta-u^*)
&=h_X^d\|r_\Omega(v_\theta)\|_{\ell^2(X)}^2
+h_Y^{d-4}
\|r_{\partial\Omega}(v_\theta)\|_{\ell^2(Y)}^2
=J(v_\theta-u^*).
\end{align}

\item $\gamma=0$ (standard PINN loss) degrades to unweighted collocation scaled by $h_Y^{d-4}$:
\begin{align}
J_0(v_\theta-u^*)
&=h_Y^{d-4}
\left(
\|r_\Omega(v_\theta)\|_{\ell^2(X)}^2
+\|r_{\partial\Omega}(v_\theta)\|_{\ell^2(Y)}^2
\right).
\end{align}
\end{enumerate}

% Using the objective $J_\gamma$, we define the physics-balanced
% constrained (PBC) and physics-balanced unconstrained (PBU) solutions by
% \begin{align}
% u_\theta^{\mathrm{PBC}_\gamma}
% &\in\argmin_{v_\theta\in\NNa}
% J_\gamma(v_\theta-u^*) \qquad \text{and}\qquad
% \label{eq:pbc-minimization-problem}
% u_\theta^{\mathrm{PBU}_\gamma}
% \in\argmin_{v_\theta\in\NN}
% J_\gamma(v_\theta-u^*),
% % \label{eq:pb-minimization-problem}
% \end{align}
% whenever these minimizers exist. Here, PBC refers to minimization over a stable
% subclass (e.g., with a BII constraint), whereas PBU
% refers to minimization over the full, unconstrained trial class.

\subsection{Residual-BII Ratio Vector}
% Fix $\gamma\in[0,2]$. The admissible subclass $\NNa$
% is determined by fixing the BII thresholds $\mathbf C_{\mathrm{BII}}=(C_{\mathrm{BII},\Omega},C_{\mathrm{BII},\partial\Omega})$ a priori.
% The unconstrained minimizer $u_\theta^{\mathrm{PBU}_\gamma}\in\NN$
% is computed by minimizing $J_\gamma(v_\theta-u^*)$ over the full neural class,
% independently of $\mathbf C_{\mathrm{BII}}$.
% Since the BII constraints are not imposed in this minimization,
% admissibility is not guaranteed:
% \[
% u_\theta^{\mathrm{PBU}_\gamma}\notin\NNa
% \qquad\text{in general}.
% \]
% Without admissibility, the preceding stability and convergence arguments
% cannot be applied directly. Instead of fixing $\NNa$
% first, we reverse the logic. Given a sequence
% $\{u_\theta^{\mathrm{PBU}_\gamma}\}$, we seek admissible subclasses
% $\NNa$, defined by fixed BII thresholds
% $\mathbf C_{\mathrm{BII}}$ independent of $h_\theta$, that contain the
% corresponding minimizers for all sufficiently small $h_\theta$.
% \subsection{Residual-BII Ratio Vector}
For any $\gamma\in[0,2]$, the admissible subclass $\NNa$
is determined by fixing the BII thresholds
$\mathbf C_{\mathrm{BII}}
=(C_{\mathrm{BII},\Omega},C_{\mathrm{BII},\partial\Omega})$
a priori.
The unconstrained minimizer
$u_\theta^{\mathrm{PBU}_\gamma}\in\NN$
is computed by minimizing $J_\gamma(v_\theta-u^*)$
over the full trial class, independently of $\mathbf C_{\mathrm{BII}}$.
Since the BII constraints are not imposed in this minimization,
admissibility is not guaranteed, i.e., 
\[
u_\theta^{\mathrm{PBU}_\gamma}\notin\NNa
\qquad\text{in general}.
\]

Instead of fixing $\NNa$ first, we reverse the logic.
Given a sequence $\{u_\theta^{\mathrm{PBU}_\gamma}\}$,
we seek admissible subclasses $\NNa$, defined by fixed BII
thresholds $\mathbf C_{\mathrm{BII}}$ independent of $h_\theta$,
that contain the corresponding minimizers for all sufficiently
small $h_\theta$.
For any trial function $v_\theta\in\NN$ with $r_D(v_\theta)\in H^{\beta_D}(D)$ for both $D\in\{\Omega,\partial\Omega\}$, define its BII ratio vector in $\mathbb{R}^2$ by
\begin{equation}\label{eq:residual-bii-vector}
\begin{aligned}
\mathbf B(v_\theta):=\bigl(B_\Omega(v_\theta),B_{\partial\Omega}(v_\theta)\bigr), \qquad
B_D(v_\theta):=h_\theta^{\varepsilon_*}
\frac{\|r_D(v_\theta)\|_{H^{\beta_D}(D)}}
{\|r_D(v_\theta)\|_{H^{\beta_D-\varepsilon_*}(D)}},
\end{aligned}
\end{equation}
with $B_D(v_\theta):=0$ if $r_D(v_\theta)\equiv0$ for $D\in\{\Omega,\partial\Omega\}$.

\begin{lemma}[A posteriori existence of an admissible subclass]
\label{lem:bounded-bii-admissibility}
Let $\{u_\theta^{\PB_\gamma}\}_{h_\theta}$ be a sequence of PBU minimizers.
Suppose that their BII ratio vectors $\mathbf B(u_\theta^{\PB_\gamma})$ are eventually bounded as
$h_\theta\to0$, then there exist BII constants
$\mathbf C_{\mathrm{BII}}=(C_{\mathrm{BII},\Omega},C_{\mathrm{BII},\partial\Omega})$ such that
\begin{equation}\label{eq:pbu-eventual-admissibility}
\mathbf B(u_\theta^{\PB_\gamma})\le\mathbf C_{\mathrm{BII}}
\quad\,\,\text{and hence}\quad\,\,
u_\theta^{\PB_\gamma}\in\NNa \,\text{ for all sufficiently small $h_\theta$.}
\end{equation}
\end{lemma}
\begin{proof}
Let 
$\mathcal B_R(0):=\{\mathbf b\in\mathbb R^2:\|\mathbf b\|<R\}$
for $R>0$ with any norm on $\mathbb R^2$.
Since the BII ratio vectors $\mathbf B(u_\theta^{\PB_\gamma})$ are eventually bounded, we have
\begin{equation}
\limsup_{h_\theta\to0}\|\mathbf B(u_\theta^{\PB_\gamma})\|<\infty\quad
\Longleftrightarrow \quad
\begin{gathered}
\exists \, R>0:\ \mathbf B(u_\theta^{\PB_\gamma})\in\mathcal B_R(0)\\
\text{for all sufficiently small }h_\theta.
\end{gathered}
\end{equation}
% By the hypothesis, choose such an $R$, independent of $h_\theta$.
% Equivalence of norms on  $\mathbb R^2$ gives
% $\|\mathbf B(u_\theta^{\PB_\gamma})\|_\infty
% \le c\|\mathbf B(u_\theta^{\PB_\gamma})\|$
% with $c>0$ independent of $h_\theta$.
% Set $C_{\mathrm{BII},D}=cR$ for $D\in\{\Omega,\partial\Omega\}$.
% Then the componentwise bounds hold for all sufficiently small $h_\theta$,
% and Definition~\ref{ass:endpoint-bii} yields
% \eqref{eq:pbu-eventual-admissibility}.
For such an $R$, equivalence of norms on $\mathbb R^2$ gives
$\|\mathbf B(u_\theta^{\PB_\gamma})\|_\infty
\le c\|\mathbf B(u_\theta^{\PB_\gamma})\|<cR$
for all sufficiently small $h_\theta$, with $c>0$ independent of
$h_\theta$.
Taking $C_{\mathrm{BII},D}=cR$ for $D\in\{\Omega,\partial\Omega\}$
yields \eqref{eq:pbu-eventual-admissibility} by
Definition~\ref{ass:endpoint-bii}.
\end{proof}

Once essential boundedness ensures
$u_\theta^{\PB_\gamma}\in\NNa$ a posteriori for some threshold
$\mathbf C_{\mathrm{BII}}$, the same discrete stability estimate
in Theorem~\ref{cor:physics-balanced-stability} then applies
under the corresponding sampling conditions.

% Once essential boundedness ensures $u_\theta^{\PB_\gamma} \in \NNa$ a posteriori for some threshold $\mathbf{C}_{\mathrm{BII}}$, the same discrete stability estimate in Theorem~\ref{cor:physics-balanced-stability} applies under the corresponding sampling conditions.

% With stability established on the admissible subclass $\NNa$, the remaining question is
% whether this subclass contains approximants $A_hu^*$ satisfying
% $J_\gamma(A_hu^*-u^*)\to0$ as $h_\theta\to0$.
% The next section addresses this consistency requirement and its
% implications for convergence.

\section{\texorpdfstring{Consistency and convergence}{Consistency and convergence}}
\label{sec:approximation-consistency}

With stability established on the admissible subclass $\NNa$,
the remaining question is whether this subclass contains
approximants $A_hu^*$ satisfying
$J_\gamma(A_hu^*-u^*)\to0$ as $h_\theta\to0$.

This section therefore addresses this consistency requirement
and its implications for convergence.
We first establish conditional convergence of PBC and PBU
on general trial classes, then verify the PBC approximation
hypothesis by constructing admissible comparison networks
for the prescribed tanh-MLP family.

\subsection{\texorpdfstring{Convergence of PBC and PBU on general trial classes}{Convergence of PBC and PBU on general trial classes}}
\label{subsec:admissible-comparison-criterion}
% Consider the neural class $\NN$ and its admissible
% subclass $\NNa$, with design scale $h_\theta$.
% The BII orders and thresholds are fixed independently of $h_\theta$.
In this subsection, sufficient oversampling matched to $\mathbf C_{\mathrm{BII}}$
means that $X$ and $Y$ are quasi-uniform with uniformly bounded
mesh ratios, satisfy
\begin{equation}\label{eq:sampling-scale-lower-bound}
h_X\le h_Y<1,\qquad
h_\theta\le\kappa_S h_S,\quad S\in\{X,Y\},
\end{equation}
for constants $\kappa_S\ge1$ independent of $h_\theta$,
and, together with the chosen admissible class, satisfy the hypotheses of either
Lemma~\ref{thm:standard-residual-stability} or Lemma~\ref{thm:int-residual} for these thresholds.

With stability established on the admissible subclass $\NNa$, we next examine
consistency and convergence, beginning with the PBC solution.

\begin{theorem}[Conditional convergence of PBC]
\label{thm:general-bii-consistency}
Suppose Assumption~\ref{ass:pde} holds.
Let $u^*\in H^m(\Omega)$ be the solution of 
\eqref{eq:PDE} with $m>2+d/2$.
Fix the admissible class $\NNa$ with threshold
$\mathbf C_{\mathrm{BII}}=(C_{\mathrm{BII},\Omega},C_{\mathrm{BII},\partial\Omega})$.
For any loss parameter $\gamma\in[0,2]$, under sufficient oversampling
$h_X,h_Y=\mathcal O(h_\theta)$ matched to $\mathbf C_{\mathrm{BII}}$,
suppose there exists an admissible comparison approximant 
$A_hu^*\in\NNa$ satisfying
\begin{equation}\label{eq:general-bii-comparison-residuals}
\|u^*-A_hu^*\|_{W^{2,\infty}(\Omega)}\le Ch_\theta^{p_2},
\qquad
\|u^*-A_hu^*\|_{L^\infty(\Omega)}\le Ch_\theta^{p_0},
\end{equation}
with $p_2>2-\gamma$ and $p_0>3/2$.
Then the PBC minimizer satisfies
\begin{equation}\label{eq:general-bii-matched-error}
\|u_\theta^{\mathrm{PBC}_\gamma}-u^*\|_{H^2(\Omega)}^2
\lesssim J_\gamma(A_hu^*-u^*) \lesssim h_\theta^{2\min\{p_2+\gamma-2,\, p_0-3/2\}},
\end{equation}
where the constant is independent of $h_\theta,h_X,h_Y$.
% The minimizers converge to $u^*$ in $H^2(\Omega)$ as
% $h_\theta\to0$ if $p_2>2-\gamma$ and $p_0>3/2$.
\end{theorem}
\begin{proof}
The definition of $J_\gamma$ in
\eqref{eq:parameterized-discrete-stability} and the bounds in
\eqref{eq:general-bii-comparison-residuals} give
\begin{align}\label{eq:admissible-upper-sampling}
J_\gamma&(A_hu^*-u^*)\\
&\lesssim h_X^{\frac d2\gamma}
h_Y^{d-4+2\gamma-\frac d2\gamma}|X|\,
\|\mathcal L(A_hu^*-u^*)\|_{L^\infty(\Omega)}^2  +h_Y^{d-4}|Y|\,
\|A_hu^*-u^*\|_{L^\infty(\partial\Omega)}^2\notag\\
&\lesssim h_X^{\frac d2\gamma}
h_Y^{d-4+2\gamma-\frac d2\gamma}|X|\,
\|A_hu^*-u^*\|_{W^{2,\infty}(\Omega)}^2 +h_Y^{d-4}|Y|\,
\|A_hu^*-u^*\|_{L^\infty(\Omega)}^2\notag\\
&\lesssim h_X^{\frac d2\gamma-d}h_Y^{d-4+2\gamma-\frac d2\gamma}h_\theta^{2p_2}
+h_Y^{-3}h_\theta^{2p_0}\notag\\
&\lesssim h_\theta^{2(p_2+\gamma-2)}
+h_\theta^{2(p_0-3/2)}.
 \notag
\end{align}
Sobolev embedding gives $u^*\in W^{2,\infty}(\Omega)$ and
quasi-uniformity gives $|X|\approx h_X^{-d}$ and
$|Y|\approx h_Y^{-(d-1)}$, where $|X|$ and $|Y|$ denote the numbers of
interior and boundary collocation points, respectively. The last inequality follows from
\eqref{eq:sampling-scale-lower-bound}.

Since $u_\theta^{\mathrm{PBC}_\gamma}\in\NNa$ and
$A_hu^*\in\NNa$,
combining \eqref{eq:parameterized-discrete-stability},
\eqref{eq:pbc-minimization-problem} and
\eqref{eq:admissible-upper-sampling} yields
\begin{align}\label{eq:general-bii-error-bound}
\|u_\theta^{\mathrm{PBC}_\gamma}-u^*\|_{H^2(\Omega)}^2
\lesssim J_\gamma(u_\theta^{\mathrm{PBC}_\gamma}-u^*)
&\lesssim J_\gamma(A_hu^*-u^*)\\
&\lesssim h_\theta^{2\min\{p_2+\gamma-2,\,p_0-3/2\}}.
\notag
\end{align}
The proof of \eqref{eq:general-bii-matched-error} is complete.
\end{proof}

% $u_\theta^{\PB_\gamma} \in \NNa$ holds a posteriori via bounded ratio vector $\mathbf{B}(u_\theta^{\PB_\gamma}) \le \mathbf{C}_{\mathrm{BII}}$.
% Different form PBC, a comparison network $A_hu^*\in\NN$ suffices, since PBU
% minimizes over the full neural class. 
% The following theorem gives the
% corresponding conditional convergence result.

For the PBC problem, consistency requires an accurate comparison approximant
$A_hu^*\in\NNa$.
For the prescribed $\tanh$ family,
Subsection~\ref{subsec:bii-consistency} verifies this approximation
hypothesis through a localized correction.

For the PBU problem, a comparison function $A_hu^*\in\NN$ suffices,
since PBU minimizes over the full trial class.
However, admissibility is not guaranteed, so the preceding
stability and convergence arguments cannot be applied directly.
Instead of fixing $\NNa$ first, we reverse the logic.
Under essential boundedness of the BII ratio vectors,
Lemma~\ref{lem:bounded-bii-admissibility} provides fixed thresholds
$\mathbf C_{\mathrm{BII}}$ independent of $h_\theta$.
Then, for all sufficiently small $h_\theta$,
$u_\theta^{\PB_\gamma}\in\NNa$ holds a posteriori via
the bounded ratio vector
$$\mathbf B(u_\theta^{\PB_\gamma})\le\mathbf C_{\mathrm{BII}}.$$
Under the corresponding sampling conditions, the same discrete
stability estimate in Theorem~\ref{cor:physics-balanced-stability}
applies.
The following theorem gives the corresponding conditional
convergence result.

\begin{theorem}[Conditional convergence of PBU]
\label{thm:pb-conditional-convergence}
Let $R>0$. For any $\gamma\in[0,2]$, assume the empirical sequence  $\{u_\theta^{\PB_\gamma}\}_{h_\theta}$ satisfies $\mathbf{B}(u_\theta^{\PB_\gamma}) \in \mathcal{B}_R(0)$ essentially, fixing an admissible class $\NNa$ with $\mathbf{C}_{\mathrm{BII}}$ matching $R$.  If there exists a comparison approximant $A_hu^*\in\NN$ satisfying $W^{2,\infty}$ and $L^\infty$ bounds \eqref{eq:general-bii-comparison-residuals} with $p_2 > 2-\gamma$ and $p_0 > 3/2$, then
\begin{align}
\|u_\theta^{\mathrm{PBU}_\gamma}-u^*\|_{H^2(\Omega)}^2
\lesssim  J_\gamma(A_hu^*-u^*) \lesssim h_\theta^{2\,\min\{p_2+\gamma-2,\;p_0-3/2\}}.
\end{align}
\end{theorem}

\begin{proof}
By Lemma~\ref{lem:bounded-bii-admissibility}, the PBU minimizers
belong to $\NNa$ for all sufficiently small $h_\theta$.
The stability estimate \eqref{eq:parameterized-discrete-stability},
global minimality in \eqref{eq:pbc-minimization-problem},
and \eqref{eq:admissible-upper-sampling} therefore yield
\begin{align*}
\|u_\theta^{\mathrm{PBU}_\gamma}-u^*\|_{H^2(\Omega)}^2
\lesssim J_\gamma(u_\theta^{\mathrm{PBU}_\gamma}-u^*) &\le J_\gamma(A_hu^*-u^*)\\
&\lesssim h_\theta^{2\min\{p_2+\gamma-2,\,p_0-3/2\}}.
\end{align*}
\end{proof}

\subsection{\texorpdfstring{Convergence based on tanh-MLPs}{Convergence  based on tanh-MLPs}}
\label{subsec:neural-approximation-power}
\label{subsec:bii-consistency}
In this subsection, we specialize $\NN$ to fixed-depth,
fully connected $\tanh$ networks with $L-1\ge1$ hidden layers.
The same residual-BII conditions and thresholds define $\NNa$.
For each prescribed architecture, let $n_\ell$ denote the width of
layer $\ell$, with $n_0=d$ and $n_L=1$.
The network is given by
\begin{equation}\label{eq:network-architecture}
% \begin{aligned}
v_\theta=W_Lz_{L-1}+b_L \quad \text{with}\quad  z_\ell=\sigma_\ell(W_\ell z_{\ell-1}+b_\ell) \text{ for } 1\le\ell<L \text{ and }z_0=x,
% \end{aligned}
\end{equation}
where $\sigma_\ell=\tanh$ acts componentwise.
We define the design scale in terms of the parameter count $n_P=\sum_{\ell=1}^{L}(n_{\ell-1}+1)n_\ell$ by
\begin{equation}\label{eq:number_hidden_units}
h_\theta=c_\sigma n_P^{-1/d},
% \qquad
% n_P=\sum_{\ell=1}^{L}(n_{\ell-1}+1)n_\ell,
\end{equation}
where $c_\sigma>0$ is fixed and $n_P$ counts all weight and bias
positions, including entries equal to zero.
The number of hidden neurons is $n_H=\sum_{\ell=1}^{L-1}n_\ell$, and
$n_{P_0}(\theta)=\|\theta\|_0\le n_P$ counts the nonzero parameters.

% Extend $u^*$ to a containing
% cube using \cite[Theorem~5.24]{AdamsFournier2003} and
% \cite[Theorem~1.122(ii)]{Triebel2006FunctionSpacesIII}.
% Adapting \cite[Theorem~1]{Park2026} to this architecture using
% shared polynomial features, with total parameter count comparable
% to the polynomial-space dimension. 
% For $u^*\in H^m(\Omega)$ with $m\ge4$, the $p=2$ case provides
% an accurate network $\widetilde u^*\in\mathcal N_\theta$ such that
% \begin{equation}\label{eq:simultaneous-tanh-sobolev-approximation}
% \|u^*-\widetilde u^*\|_{H^j(\Omega)}
% \lesssim h_\theta^{m-j},\qquad j\in\{0,2,4\}.
% \end{equation}
% By the interpolation inequality
% $\|v\|_{W^{k,\infty}(\Omega)}\lesssim
% \|v\|_{L^2(\Omega)}^{1-a_k}\|v\|_{H^4(\Omega)}^{a_k}$
% for $v\in H^4(\Omega)$, where $a_k=(k+d/2)/4$ and
% $k\in\{0,2\}$ \cite{LiZhang2022GN}, the preceding estimates give
% \begin{equation}\label{eq:simultaneous-tanh-approximation}
% \|u^*-\widetilde u^*\|_{W^{k, \infty}(\Omega)}
% \lesssim n_P^{-\frac{m-k-d/2}{d}} \lesssim  h_\theta^{m-k-\frac{d}{2}},\quad k\in\{0,2\}
% \end{equation}
% Related estimates in terms of $n_{P_0}$ are given in
% \cite[Proposition~4.8]{GuhringRaslan2021}.
% The base network uses at most $n_\ell-2$ neurons per hidden layer,
% reserving two for the localized correction below; since
% $n_\ell-2 \approx n_\ell$, the approximation orders remain unchanged.
% \par

For the approximation results below, we increase the common hidden-layer
width and allow unrestricted weights and biases. Park's construction
\cite{Park2026} can be extended to this family using shared polynomial
features, with $n_P \approx N^d$ total parameters for coordinatewise
polynomial degree $N$.
For $u^*\in H^m(\Omega)$ with $m\ge4$, applying this construction
after Sobolev extension to a containing cube
(see \cite[Theorem~5.24]{AdamsFournier2003} and
\cite[Theorem~1.122(ii)]{Triebel2006FunctionSpacesIII})
gives a single network $\widetilde u^*\in\NN$ satisfying the following
standard simultaneous approximation estimates
\begin{equation}\label{eq:simultaneous-tanh-approximation}
\|u^*-\widetilde u^*\|_{H^k(\Omega)}
\lesssim n_P^{-(m-k)/d} \approx h_\theta^{m-k},
\qquad k\in\{0,2,4\}.
\end{equation}
Applying the Gagliardo--Nirenberg inequality
\cite[Theorem~1.3]{LiZhang2022GN} to \eqref{eq:simultaneous-tanh-approximation} then gives
\begin{equation}\label{eq:simultaneous-tanh-maxnorm}
\|u^*-\widetilde u^*\|_{W^{k,\infty}(\Omega)}
\lesssim 
\|u^*-\widetilde u^*\|_{H^k(\Omega)}^{1-d/4}
\|u^*-\widetilde u^*\|_{H^{k+2}(\Omega)}^{d/4} \lesssim h_\theta^{m-k-d/2},\quad k\in\{0,2\}.
\end{equation}
For related estimates in terms of the nonzero parameter count $n_{P_0}$,
see \cite{GuhringRaslan2021}.

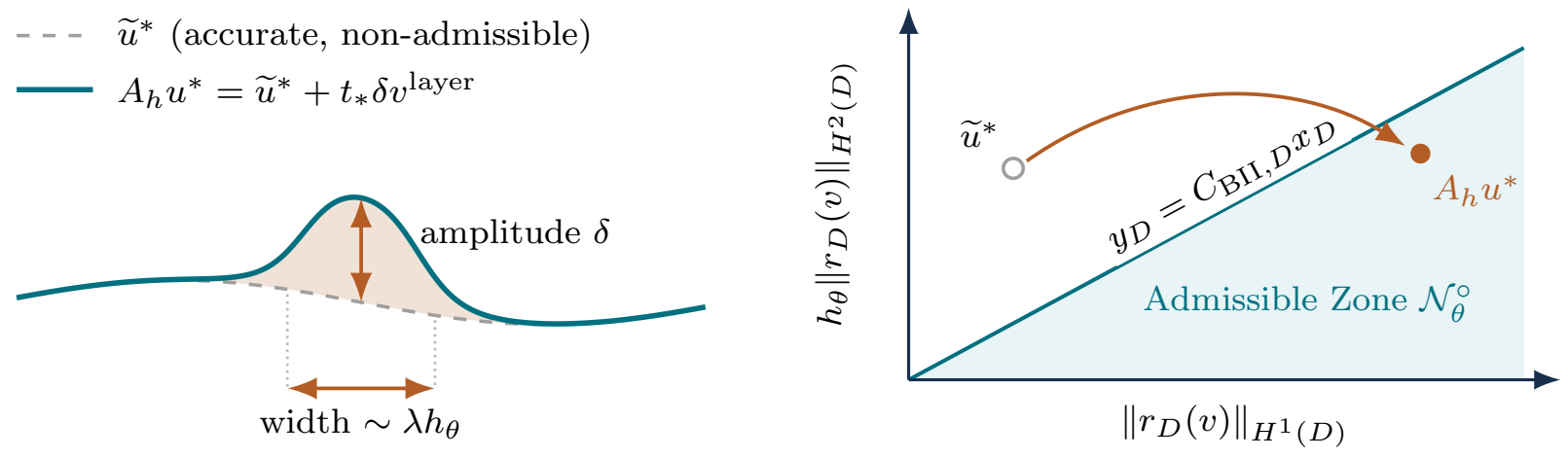
\begin{figure}[t]
\centering
\newcommand{\biiFigureScale}{0.85}
\pgfmathdeclarefunction{biitanh}{1}{%
\pgfmathparse{sign(#1)*(1-exp(-2*abs(#1)))/(1+exp(-2*abs(#1)))}}
\pgfmathdeclarefunction{biibase}{1}{%
\pgfmathparse{1.65+0.18*sin(60*#1)}}
\pgfmathdeclarefunction{biilayer}{1}{%
\pgfmathparse{0.85*(biitanh(biitanh((#1-3)/0.46+1))
-biitanh(biitanh((#1-3)/0.46-1)))
/(2*biitanh(biitanh(1)))}}

\resizebox{\biiFigureScale\linewidth}{!}{
\begin{tikzpicture}[
x=1cm,y=1cm,>=Latex,
font=\fontsize{8.5}{10}\selectfont,
base/.style={draw=gray!75,dashed,line width=0.8pt},
corrected/.style={draw=sea,line width=1.25pt},
correction/.style={draw=warm,line width=0.85pt},
title/.style={font=\fontsize{9}{11}\selectfont},
every node/.style={inner sep=2pt}]

\path[use as bounding box] (0,-0.40) rectangle (13,4.3);
% (a) Localized correction: a spatial cross-section.
\draw[base] (0.20,3.82)--(0.82,3.82);
\node[anchor=west] at (0.94,3.82) {$\widetilde u^*$ (accurate, non-admissible)};
\draw[corrected] (0.20,3.38)--(0.82,3.38);
\node[anchor=west] at (0.94,3.38) {$A_hu^*=\widetilde u^*+t_*\delta v^{\mathrm{layer}}$};

\fill[warm!18]
plot[domain=0.2:5.8,samples=181] (\x,{biibase(\x)})
-- plot[domain=5.8:0.2,samples=181]
(\x,{biibase(\x)+biilayer(\x)}) -- cycle;
\draw[base] plot[domain=0.2:5.8,samples=181] (\x,{biibase(\x)});
\draw[corrected] plot[domain=0.2:5.8,samples=181] (\x,{biibase(\x)+biilayer(\x)});

\draw[correction,<->] (3,1.65)--(3,2.50);
\node[anchor=west] at (3.4,2.2) {amplitude $\delta$};
\draw[gray!55,densely dotted,line width=0.6pt]
(2.4,{biibase(2.4)-0.04})--(2.4,0.95);
\draw[gray!55,densely dotted,line width=0.6pt]
(3.6,{biibase(3.6)-0.04})--(3.6,0.95);
\draw[correction,<->] (2.4,0.95)--(3.6,0.95);
\node[anchor=north] at (3,0.86) {width $\sim\lambda h_\theta$};
\node[title,align=center,anchor=north] at (3.05,0.38)
{(a) Preserved approximation power\\[-1pt]
{\fontsize{8}{10}\selectfont via localized correction}};

% (b) For a fixed D, the image of its BII condition in the norm plane.
% The curved arrow is schematic, not a monotonicity assertion in t.
\begin{scope}[shift={(7.45,1.02)}]
\fill[pale] (0,0)--(5,0)--(5,2.70)--cycle;
\draw[sea,line width=0.9pt] (0,0)--(5,2.70);
\draw[->,draw=ink,line width=0.7pt] (0,0)--(5.30,0);
\draw[->,draw=ink,line width=0.7pt] (0,0)--(0,3.03);

\node[anchor=north] at (2.65,-0.13)
{$\|r_D(v)\|_{H^1(D)}$};
\node[rotate=90] at (-0.59,1.52)
{$h_\theta\|r_D(v)\|_{H^2(D)}$};
\node[rotate=28.37,anchor=south,fill=white,inner sep=1.2pt]
at (2.65,1.431) {$y_D=C_{\mathrm{BII},D}x_D$};
\node[text=sea,align=center] at (3.25,0.62)
{Admissible Zone $\NNa$};

\draw[draw=gray!75,fill=white,line width=0.8pt]
(0.85,1.72) circle (2.3pt);
\node[anchor=south east] at (0.82,1.82) {$\widetilde u^*$};
\fill[warm] (4.16,1.84) circle (2.3pt);
\node[anchor=north west,text=warm] at (4.18,1.74) {$A_hu^*$};
\draw[correction,->]
(0.96,1.78) .. controls (1.95,2.48) and (3.15,2.45)
.. (4.04,1.90);
\end{scope}
\node[title,anchor=north] at (10.05,0.38)
{(b) BII admissibility};
\end{tikzpicture}
}
\caption{The dual objective of localized layer correction.
Starting from an accurate base candidate
$\widetilde u^*\in\NN$, we add a thin, localized layer
correction to obtain an admissible comparison network
$A_hu^*\in\NNa$.
(a) The localized correction keeps the boundary error
small enough not to ruin the approximation order.
(b) The correction forces the residual BII ratios
into the admissible zone.}
\label{fig:bii-admissible-correction}
\vspace{-15pt}
\end{figure}

Standard theory gives an accurate $\widetilde u^*\in\NN$,
but we need an admissible comparison network $A_hu^*\in\NNa$.
The strategy is a localized layer correction.
We take $\widetilde u^*\in\NN$ as the base candidate with two neurons in each
hidden layer reserved for the correction and,
without loss of generality, assume
$\widetilde u^*\notin\NNa$; otherwise, we simply set
$A_hu^*=\widetilde u^*$.
We add a thin, localized boundary layer correction:
\begin{align}
A_hu^*=\widetilde u^*+t_*\delta v^{\mathrm{layer}}.
\end{align}
The correction forces the residual BII ratios into the
admissible zone while keeping the boundary error small
enough not to ruin the approximation order. For computational simplicity, we present only the case
$\beta_\Omega=\beta_{\partial\Omega}=2$ with $\varepsilon_*=1$,
the lowest integer-order BII setting with an integer gap
sufficient for the stability argument in
Lemma~\ref{thm:int-residual}.
Higher orders can be treated similarly. As illustrated in
Figure~\ref{fig:bii-admissible-correction},
we seek $\lambda$, $\delta$, and $t_*\in\{1,2,3\}$
such that $A_hu^*$ achieves the following:
% \begin{enumerate}
% \renewcommand{\labelenumi}{(\roman{enumi})}
% \item Admissibility:
\begin{align*}
&\text{(i) Admissibility:} \quad  h_\theta
\frac{\|r_D(A_hu^*)\|_{H^2(D)}}
{\|r_D(A_hu^*)\|_{H^1(D)}}
\le C_{\mathrm{BII},D}
\qquad  \text{for both } D\in\{\Omega,\partial\Omega\}.\\
% \end{align*}
% % \item Preserved Approximation Power:
% \begin{align*}
&\text{(ii) Preserved approximation power:} \quad  \|u^*-A_hu^*\|_{W^{k,\infty}(\Omega)}
\approx
\|u^*-\widetilde u^*\|_{W^{k,\infty}(\Omega)}.
\end{align*}
% \end{enumerate}
For any prescribed positive BII thresholds, the following
lemma constructs such a comparison network $A_hu^*$ within
the same prescribed architecture.

\begin{lemma}[Convergence of an admissible comparison network]
\label{lem:bii-feasible-comparison}
Let $d\le3$ and assume $m\ge4$.
For any prescribed threshold $\mathbf{C}_{\mathrm{BII}}$ with integer BII
order $\beta_\Omega=\beta_{\partial\Omega}=2$ and integer  gap $\varepsilon_* = 1$, and sufficiently small $h_\theta$, the repaired approximant $A_h u^* := \widetilde{u}^* + t_*\delta v^{\mathrm{layer}}$ satisfies:
\begin{align*}
A_h u^* \in \NNa,
\end{align*}
with comparison rates $p_2$ and $p_0$:
\begin{align}\label{eq:bii-comparison-maxnorm}
\|u^*-A_hu^*\|_{W^{2,\infty}(\Omega)} & \lesssim h_\theta^{p_2},
&&p_2 := m - 2 - d/2, \\
\|u^*-A_hu^*\|_{L^\infty(\Omega)} & \lesssim  h_\theta^{p_0},  &&p_0 := m - d/2.
\end{align}
\end{lemma}

\begin{proof}
For $d\ge2$, define the candidate family
$v_t:=\widetilde u^*+\delta t\,v^{\mathrm{layer}}$ for
$t\in\{1,2,3\}$, where $v^{\mathrm{layer}}$ is a localized ridge
layer of width $\lambda h_\theta$ along direction $\tau$, given by
\begin{equation}\label{eq:correction-ridge-definition}
v^{\mathrm{layer}}(x)
=
\chi\!\left(\frac{\tau\cdot x-a}{\lambda h_\theta}\right),
\qquad
\chi(s)
=
\tanh^{\circ(L-1)}(s+1)-\tanh^{\circ(L-1)}(s-1).
\end{equation}
The layer is realized by two reserved parallel neuron chains along
a transverse slice crossing $\Omega$ and a smooth patch of
$\partial\Omega$.
Each derivative contributes $(\lambda h_\theta)^{-1}$
and $L^2$ scaling contributes $(\lambda h_\theta)^{1/2}$.
Integration across the transverse ridge gives
\begin{align}\label{eq:correction-lower-order-scaling}
\|\mathcal L v^{\mathrm{layer}}\|_{H^j(\Omega)}
&\simeq(\lambda h_\theta)^{-j-3/2},
&
\|v^{\mathrm{layer}}\|_{H^j(\partial\Omega)}
&\simeq(\lambda h_\theta)^{-j+1/2}, \quad j\in\{1,2\}.
\end{align}

\smallskip\noindent
\textit{Step 1: Write the ratios and isolate the risk.}
By linearity, the residuals satisfy
\[
r_\Omega(v_t)
=
r_\Omega(\widetilde u^*)+\delta t\,\mathcal L v^{\mathrm{layer}},
\qquad
r_{\partial\Omega}(v_t)
=
r_{\partial\Omega}(\widetilde u^*)+\delta t\,v^{\mathrm{layer}}|_{\partial\Omega}.
\]
Let $w_\Omega=\mathcal L v^{\mathrm{layer}}$ and
$w_{\partial\Omega}=v^{\mathrm{layer}}|_{\partial\Omega}$.
Applying the triangle inequality only to the numerator gives, for $D\in\{\Omega,\partial\Omega\}$, BII order
$\beta_D=2$ and gap $\varepsilon_*=1$:
\[
h_\theta
\frac{\|r_D(v_t)\|_{H^2(D)}}{\|r_D(v_t)\|_{H^1(D)}}
\le
\frac{
h_\theta\|r_D(\widetilde u^*)\|_{H^2(D)}
+\delta t\,h_\theta\|w_D\|_{H^2(D)}
}{
\|r_D(\widetilde u^*)+\delta t\,w_D\|_{H^1(D)}
}.
\]
The numerators grow at most linearly with $t$.
The danger is destructive cancellation in the denominator.
If $\delta t\,w_D$ cancels or nearly cancels
$r_D(\widetilde u^*)$, the denominator may vanish or become
arbitrarily small.

\smallskip\noindent
\textit{Step 2: Candidate selection protects both denominators.}
We seek a single candidate $t_*\in\{1,2,3\}$ satisfying both
lower bounds simultaneously:
\begin{align}
\|r_\Omega(v_{t_*})\|_{H^1(\Omega)}
&\ge
\frac{\delta}{2}
\|\mathcal L v^{\mathrm{layer}}\|_{H^1(\Omega)}, \quad
\|r_{\partial\Omega}(v_{t_*})\|_{H^1(\partial\Omega)}
\ge
\frac{\delta}{2}
\|v^{\mathrm{layer}}\|_{H^1(\partial\Omega)}.
\label{eq:correction-denominator-bounds}
\end{align}
Such a common $t_*$ always exists.
To see this, suppose two distinct candidates
$t_1,t_2\in\{1,2,3\}$ both fail the interior bound:
\begin{align*}
\|r_\Omega(v_{t_1})\|_{H^1(\Omega)}<
\frac{\delta}{2}
\|\mathcal L v^{\mathrm{layer}}\|_{H^1(\Omega)}, \quad
\|r_\Omega(v_{t_2})\|_{H^1(\Omega)}
<\frac{\delta}{2}
\|\mathcal L v^{\mathrm{layer}}\|_{H^1(\Omega)}.
\end{align*}
Their difference is exactly
$r_\Omega(v_{t_1})-r_\Omega(v_{t_2})
=
\delta(t_1-t_2)\mathcal L v^{\mathrm{layer}}$.
Since $|t_1-t_2|\ge1$, the triangle inequality gives
\begin{align*}
\delta\|\mathcal L v^{\mathrm{layer}}\|_{H^1(\Omega)}
\le
\|r_\Omega(v_{t_1})-r_\Omega(v_{t_2})\|_{H^1(\Omega)}
&\le
\|r_\Omega(v_{t_1})\|_{H^1(\Omega)}
+
\|r_\Omega(v_{t_2})\|_{H^1(\Omega)}\\
&<
\delta\|\mathcal L v^{\mathrm{layer}}\|_{H^1(\Omega)},
\end{align*}
which gives a contradiction.
Thus at most one candidate $t\in\{1,2,3\}$ can fail the bound on $\Omega$,
and the same argument applies on $\partial\Omega$.
By the pigeonhole principle, at least one of the three candidates
satisfies both lower bounds.

\smallskip\noindent
\textit{Step 3: Substitute the denominator bounds into the ratios.}
By \eqref{eq:correction-lower-order-scaling}, there is a constant
$K>0$, independent of $\lambda$ and $h_\theta$, such that,
for each fixed $\lambda\ge1$ and sufficiently small $h_\theta$,
$\lambda h_\theta\|w_D\|_{H^2(D)}/\|w_D\|_{H^1(D)}\le K$
for $D\in\{\Omega,\partial\Omega\}$.
The triangle inequality and
\eqref{eq:correction-denominator-bounds}, together with $t_*\le3$,
give, for $D\in\{\Omega,\partial\Omega\}$,
\begin{equation}\label{eq:repaired-residual-bii-ratio}
\begin{aligned}
h_\theta
\frac{\|r_D(v_{t_*})\|_{H^2(D)}}
{\|r_D(v_{t_*})\|_{H^1(D)}}
&\le
\underbrace{
\frac{2h_\theta}{\delta}
\frac{\|r_D(\widetilde u^*)\|_{H^2(D)}}
{\|w_D\|_{H^1(D)}}
}_{\text{Term (I): controlled by }\delta}
\quad+
\underbrace{
6h_\theta
\frac{\|w_D\|_{H^2(D)}}
{\|w_D\|_{H^1(D)}}
}_{\text{Term (II): controlled by }6K/\lambda}.
\end{aligned}
\end{equation}
Term (II) dependson the ridge layer construction and is
independent of $\widetilde u^*$ and $\delta$.

\smallskip\noindent
\textit{Step 4: Fix parameters $\lambda\to\delta\to t_*$.}
Let $C_{\mathrm{BII}}
=
\min\{C_{\mathrm{BII},\Omega},
      C_{\mathrm{BII},\partial\Omega}\}>0$.
First fix the width parameter $\lambda$ by choosing
\[
\lambda\ge
\max\left\{1,\frac{8K}{C_{\mathrm{BII}}}\right\},
\]
and take $h_\theta$ sufficiently small that
 term (II) is bounded by $3C_{\mathrm{BII}}/4$.

Next fix the amplitude $\delta$ by setting
\begin{equation}\label{eq:correction-parameter-choice}
\delta
=
\frac{8h_\theta}{C_{\mathrm{BII}}}
\max\left\{
\frac{\|r_\Omega(\widetilde u^*)\|_{H^2(\Omega)}}
{\|\mathcal Lv^{\mathrm{layer}}\|_{H^1(\Omega)}},
\frac{\|r_{\partial\Omega}(\widetilde u^*)\|_{H^2(\partial\Omega)}}
{\|v^{\mathrm{layer}}\|_{H^1(\partial\Omega)}}
\right\}.
\end{equation}
Then Term (I) is bounded by $C_{\mathrm{BII}}/4$.
Finally, choose $t_*$ as in Step 2 for this $\delta$ and set
$
A_hu^*
:=
v_{t_*}
=
\widetilde u^*+t_*\delta v^{\mathrm{layer}}$.
The two term bounds yield
\[
\mathbf B(A_hu^*)\le\mathbf C_{\mathrm{BII}}
\quad\Longleftrightarrow\quad
A_hu^*\in\NNa.
\]

\smallskip\noindent
\textit{Step 5: Preservation of approximation power in $W^{k,\infty}$.}
For $k\in\{0,2\}$, the triangle inequality gives
\begin{equation}\label{Wkinf}
\|u^*-A_hu^*\|_{W^{k,\infty}(\Omega)}
\le
C\|u^*-\widetilde u^*\|_{W^{k,\infty}(\Omega)}
+
t_*\delta\|v^{\mathrm{layer}}\|_{W^{k,\infty}(\Omega)}.
\end{equation}
The first term is bounded by \eqref{eq:simultaneous-tanh-maxnorm}.
For the correction term, recall that $t_*\le3$.
The amplitude choice and layer scaling, together with the trace bound
$\|r_{\partial\Omega}(\widetilde u^*)\|_{H^2(\partial\Omega)}
\le C\|u^*-\widetilde u^*\|_{H^{5/2}(\Omega)}$,
\eqref{eq:simultaneous-tanh-approximation}, and Sobolev interpolation,
give
\begin{equation}\label{eq:correction-amplitude-bound}
\begin{aligned}
\delta
&\lesssim
h_\theta\left(
h_\theta^{5/2}
\|r_\Omega(\widetilde u^*)\|_{H^2(\Omega)}
+
h_\theta^{1/2}
\|r_{\partial\Omega}(\widetilde u^*)\|_{H^2(\partial\Omega)}
\right)
\\
&\lesssim 
h_\theta^{7/2}
\|u^*-\widetilde u^*\|_{H^4(\Omega)}
+
h_\theta^{3/2}
\|u^*-\widetilde u^*\|_{H^{5/2}(\Omega)}
\lesssim h_\theta^{m-1}.
\end{aligned}
\end{equation}
Each derivative of the ridge layer contributes
$(\lambda h_\theta)^{-1}$.
Since $\lambda$ is fixed, this gives
$\|v^{\mathrm{layer}}\|_{W^{k,\infty}(\Omega)}
\le Ch_\theta^{-k}$ for $k\in\{0,2\}$.
Combining these estimates with \eqref{Wkinf} gives
\begin{equation}
\|u^*-A_hu^*\|_{W^{k,\infty}(\Omega)}
\le
C\left(
h_\theta^{m-k-d/2}
+
h_\theta^{m-1-k}
\right)\le Ch_\theta^{m-k-d/2},
\quad k\in\{0,2\}.
\end{equation}
The last inequality uses $d/2\ge1$, proving
\eqref{eq:bii-comparison-maxnorm}.
For $d=1$, the boundary BII is automatic for sufficiently small
$h_\theta$, while an interior layer using only the interior term
in \eqref{eq:correction-parameter-choice} gives
$\delta\lesssim h_\theta^{m-1/2}$ and preserves the stated
approximation orders.
The proof is complete.

% For $d=1$, boundary BII is automatic in the endpoint $\ell^2$-norm
% for small $h_\theta$, and the same construction with an interior layer
% and only the interior term in \eqref{eq:correction-parameter-choice}
% gives $\delta\le Ch_\theta^{m-1/2}$.
% The same Gagliardo--Nirenberg estimate then yields
% \eqref{eq:bii-comparison-maxnorm}.
% The two reserved chains realize $t_*\delta v^{\mathrm{layer}}$
% with first-layer inputs $(\tau\cdot x-a)/(\lambda h_\theta)\pm1$
% and output weights $\pm t_*\delta$. Setting the connections between
% the base network and the chains to zero gives $A_hu^*\in\mathcal N_\theta$.

\end{proof}

\begin{theorem}[Convergence rates for PINNs based on $\tanh$-MLPs]
\label{thm:discrete}
Let $d\le3$, $\gamma\in[0,2]$,
and $u^*\in H^m(\Omega)$ with $m\ge4$ and $m>4+d/2-\gamma$.
For $h_\theta\simeq n_P^{-1/d}$, under sufficient
oversampling matched to the corresponding BII thresholds,
global empirical minimizers satisfy
\begin{equation}\label{eq:standard-pinn-error-decomposition}
\|u_\theta^{\mathrm{A}_\gamma}-u^*\|_{H^2(\Omega)}
\lesssim h_\theta^{m-4-d/2+\gamma}
\simeq n_P^{-\frac{m-4-d/2+\gamma}{d}},
\qquad \mathrm{A}\in\{\mathrm{PBC},\,\mathrm{PBU}\},
\end{equation}
with the following distinction:
\par\smallskip
\indent\textnormal{\textbf{PBC:}}
For any fixed positive vector of BII thresholds,
the construction yields an admissible comparison
$A_hu^*\in\NNa$ with the approximation estimates needed
for the stated rate.

\par\smallskip
\indent\textnormal{\textbf{PBU:}}
The rate holds provided the minimizer sequence
has uniformly bounded scaled residual-BII ratios, with fixed
admissibility thresholds chosen above their eventual bounds.
\end{theorem}

\begin{proof}
Lemma~\ref{lem:bii-feasible-comparison} provides
$A_hu^*\in\NNa$ with comparison rates
$p_2=m-2-d/2$ and
$p_0=m-d/2$.
Since $\gamma\in[0,2]$ and $m>4+d/2-\gamma$,
\[
\min\{p_2+\gamma-2,p_0-3/2\}
=m-4-d/2+\gamma>0.
\]
For PBC, apply Theorem~\ref{thm:general-bii-consistency}.
For PBU, the assumed scaled residual-BII ratio bounds
allow us to apply
Theorem~\ref{thm:pb-conditional-convergence}
with the same comparison network $A_hu^*\in\NN$.
Taking square roots in the resulting estimates and using
\eqref{eq:number_hidden_units} proves
\eqref{eq:standard-pinn-error-decomposition}.
\end{proof}

\begin{remark}
In \eqref{eq:standard-pinn-error-decomposition}, increasing the loss parameter $\gamma$ from $0$ to $2$ adds two powers of $h_\theta$ to the convergence exponent, corresponding to an additional factor
$h_\theta^2\simeq n_P^{-2/d}$ in the $H^2$-error bound.
The condition for a positive convergence exponent is relaxed from
$m>4+d/2$ to $m>2+d/2$, while the construction still requires $m\ge4$.

An extension to other smooth activations $\sigma$, such as GELU or
SiLU/Swish, would require quantitative simultaneous Sobolev
approximation estimates and a localized correction realizable
within the prescribed $\sigma$-network family, satisfying the
residual-BII and approximation bounds used in the proof.
\end{remark}

\section{Numerical experiments}
\label{sec:numerical-experiments}
We first describe the empirical BII implementation and common settings,
then examine smooth-solution convergence and BII control,
convergence for finite-regularity solutions, and a nonlinear extension.
We penalize empirical violations of Definition~\ref{ass:endpoint-bii} with
$\varepsilon_*=1$ on $D\in\{\Omega,\partial\Omega\}$.
For a finite weighted quadrature rule $\mathcal Q_D$ and integers $k\ge0$, set
\begin{align}\label{eq:empirical-sobolev-norm}
\|r_D(v_\theta)\|_{\widehat H^k(D)}^2
:=
\sum_{|\alpha|\le k}\sum_{(z,\omega)\in\mathcal Q_D}
\omega\left|\partial_D^\alpha r_D(v_\theta)(z)\right|^2.
\end{align}
For $\epsilon_{\rm stab}>0$, define the stabilized ratio $R_D$ and the quantity $z_D$ by
\begin{align}\label{eq:implemented-bii-fraction}
R_D(v_\theta)&:=
\frac{\|r_D(v_\theta)\|_{\widehat H^{\beta_D}(D)}}
{\|r_D(v_\theta)\|_{\widehat H^{\beta_D-1}(D)}+\epsilon_{\rm stab}},
\qquad
z_D(v_\theta):=\frac{h_\theta R_D(v_\theta)}{C_{\mathrm{BII},D}}-1.
\end{align}
Use $\mathcal B_D^+(v_\theta):=[z_D(v_\theta)]_+^2$, with
$[a]_+:=\max\{a,0\}$, or its $C^2$ smoothing
$\mathcal B_D^\delta(v_\theta):=\varphi_\delta(z_D(v_\theta))$ for $\delta>0$.
Here $\varphi_\delta(z)$ equals $0$ for $z\le0$,
$\delta^2(3t^3-3t^4+t^5)$ with $t=z/\delta$ for $0<z<\delta$,
and $z^2$ for $z\ge\delta$.
The symbol $\mathcal B_D$ denotes the selected penalty.

\begin{figure}[t]
\centering
\includegraphics[width=0.96\textwidth]{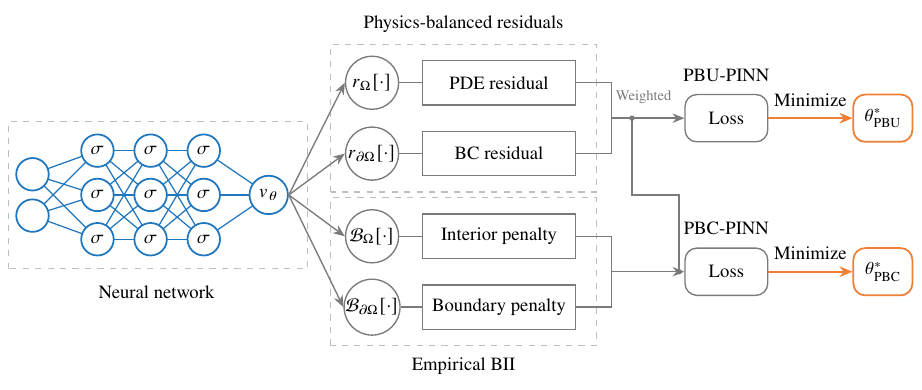}
\caption{Training objectives for PBU-PINN and PBC-PINN.}
\label{fig:bii-wls-schematic}
\vspace{-15pt}
\end{figure}

For the theoretical boundary weight $W_\gamma$, the weighted residual loss is
\begin{equation}\label{eq:numerical-residual-loss}
\mathcal J_\gamma(v_\theta):=
\|r_\Omega(v_\theta)\|_{\ell^2(X)}^2
+W_\gamma\|r_{\partial\Omega}(v_\theta)\|_{\ell^2(Y)}^2.
\end{equation}
For the linear problem \eqref{eq:PDE},
$J_\gamma(v_\theta-u^*)
=h_Y^{d-4}W_\gamma^{-1}\mathcal J_\gamma(v_\theta)$.
The training objectives have the following form, with weights scaled
according to the sampling grid.
\begin{align}\label{eq:total-training-objective}
\mathcal J_{\rm train}(v_\theta)
:=
\begin{cases}
\mathcal J_\gamma(v_\theta), & \text{PBU-PINN},\\[2pt]
\mathcal J_\gamma(v_\theta)
+\lambda_{\mathrm{BII},\Omega}\mathcal B_\Omega(v_\theta)
+\lambda_{\mathrm{BII},\partial\Omega}\mathcal B_{\partial\Omega}(v_\theta),
& \text{PBC-PINN}.
\end{cases}
\end{align}
Both methods minimize over the full trial class.
PBC-PINN adds penalties for empirical BII violations.
Figure~\ref{fig:bii-wls-schematic} summarizes this loss construction.

All numerical experiments use fully connected architectures and
the parameter-based design scale $h_\theta=n_P^{-1/2}$, where $n_P$ counts trainable
parameters. The networks have $H$ hidden layers.
The error-plot legends report
observed orders $p$ fitted to mean errors.
For PBC-PINN, we use a shared BII threshold $C_{\mathrm{BII}}$.
We set $\lambda_{\mathrm{BII},\Omega}=\lambda_{\mathrm{BII},\partial\Omega}=1$
before overall loss normalization.
Errors use composite Gauss--Legendre quadrature.
Training uses Adam followed by L-BFGS run in blocks of iterations.

\subsection{\texorpdfstring{Example 1. Smooth-solution convergence and BII control}{Example 1. Smooth-solution convergence and BII control}}
\begin{figure}[t]
\centering
\begin{subfigure}[t]{0.42\linewidth}
\centering
\includegraphics[width=\linewidth]{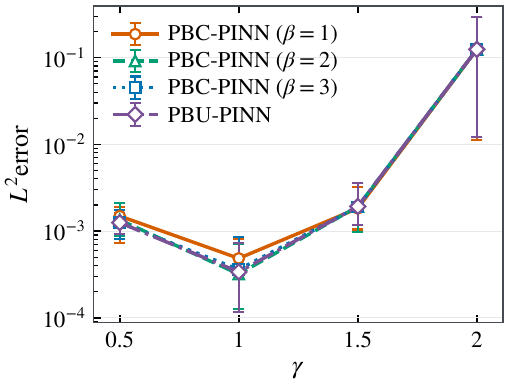}
\end{subfigure}\hspace{0.04\linewidth}%
\begin{subfigure}[t]{0.42\linewidth}
\centering
\includegraphics[width=\linewidth]{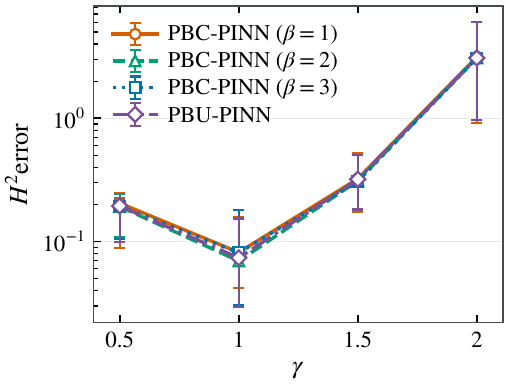}
\end{subfigure}
\caption{Smooth Poisson problem: $L^2$ errors (left) and
$H^2$ errors (right) versus the loss parameter $\gamma$. PBU-PINN
is compared with PBC-PINN at residual BII orders $\beta=1,2,3$ and threshold
$C_{\mathrm{BII}}=0.5$. All runs use two hidden layers of width 25, with
$h_X/h_\theta\approx h_Y/h_\theta\approx0.7$. Curves show means over five seeds, and error
bars span the seedwise minimum and maximum.}
\label{fig:gamma-screening}
\vspace{-15pt}
\end{figure}

On $\Omega=(0,1)^2$, consider
\begin{align}\label{eq:model-one}
-\Delta u=f \quad \text{in }\Omega \qquad 
u=0  \quad \text{on }\partial\Omega,
\end{align}
with exact solution $u^*(x,y)=\sin(\pi x)\sin(\pi y)$ and
$f(x,y)=2\pi^2\sin(\pi x)\sin(\pi y)$.
We set the residual BII orders $\beta_\Omega=\beta_{\partial\Omega}=\beta$.
We first examine how the loss parameter and residual BII
order affect accuracy, using $|X|=729$ interior midpoints and $|Y|=80$ boundary points.
Figure~\ref{fig:gamma-screening} shows that both methods have their lowest
mean errors at loss parameter $\gamma=1$.
We also test $\beta=1$ as a lower-order BII comparison, although it lies
outside the assumptions of Theorem~\ref{cor:physics-balanced-stability}.
Since $\beta=2$ is the lowest integer BII order covered by the theorem and gives
errors comparable to $\beta=3$, we adopt $(\gamma,\beta)=(1,2)$ in subsequent
experiments and denote this choice by PBC-PINN.

\begin{figure}[t]
\centering
\begin{subfigure}[t]{0.45\linewidth}
\centering
\includegraphics[width=\linewidth]{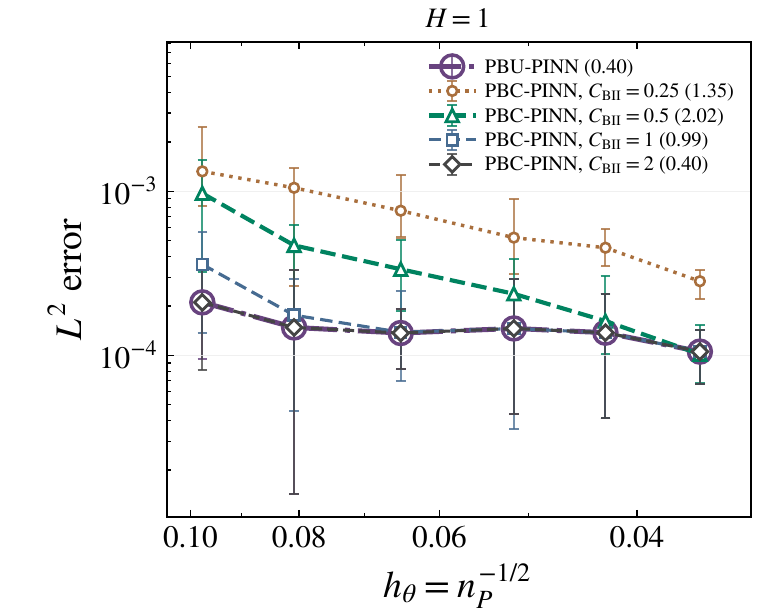}
\end{subfigure}\hspace{0.04\linewidth}%
\begin{subfigure}[t]{0.45\linewidth}
\centering
\includegraphics[width=\linewidth]{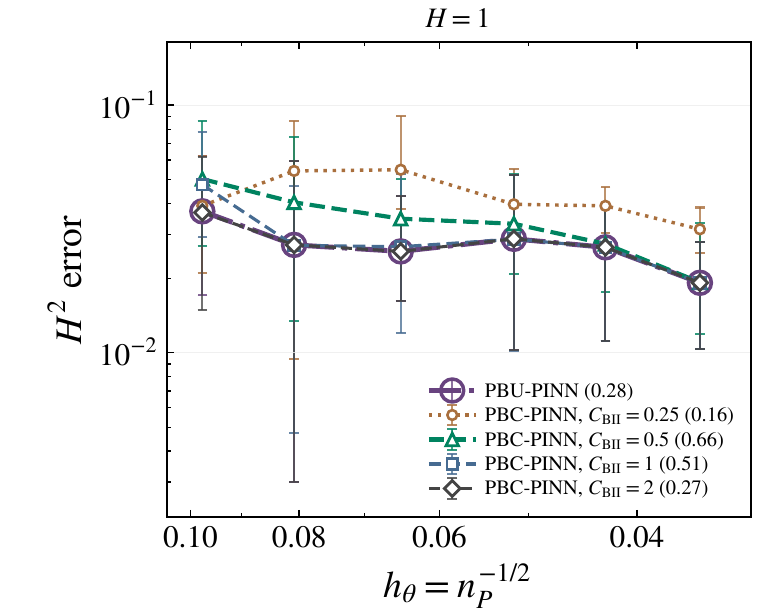}
\end{subfigure}
\par\smallskip
\begin{subfigure}[t]{0.45\linewidth}
\centering
\includegraphics[width=\linewidth]{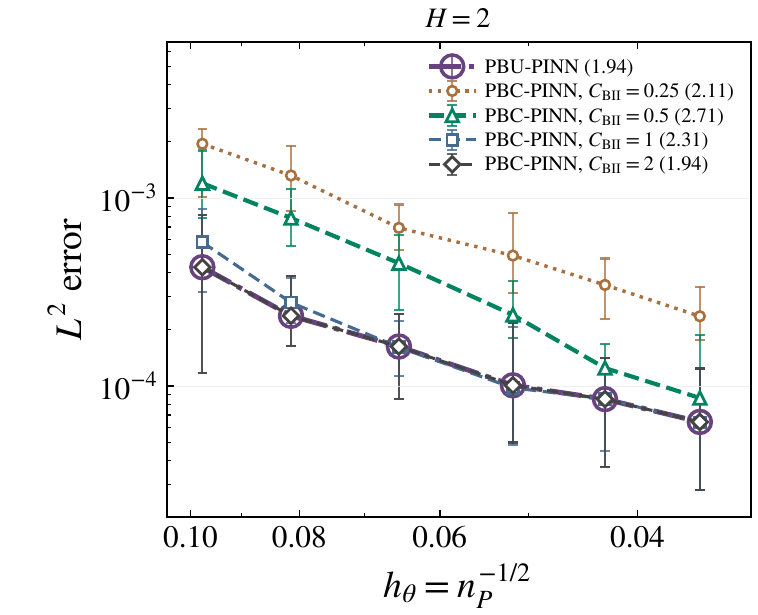}
\end{subfigure}\hspace{0.04\linewidth}%
\begin{subfigure}[t]{0.45\linewidth}
\centering
\includegraphics[width=\linewidth]{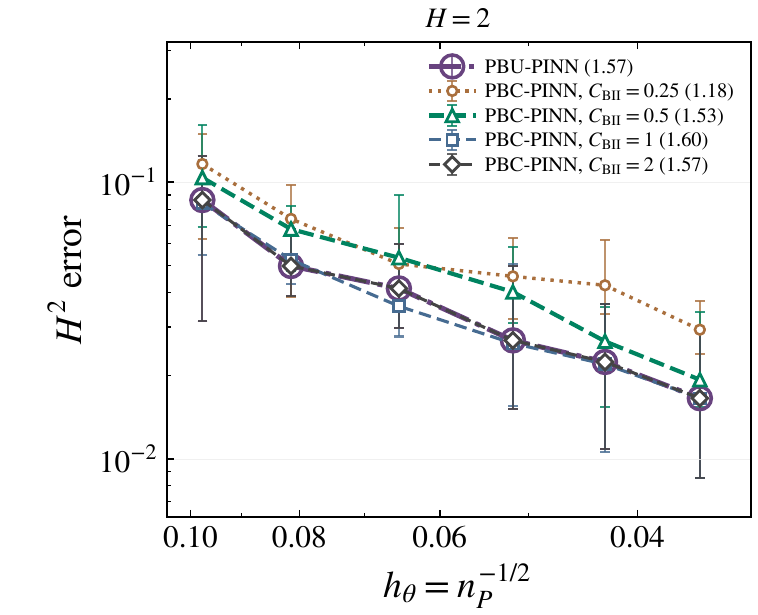}
\end{subfigure}
\par\smallskip
\begin{subfigure}[t]{0.45\linewidth}
\centering
\includegraphics[width=\linewidth]{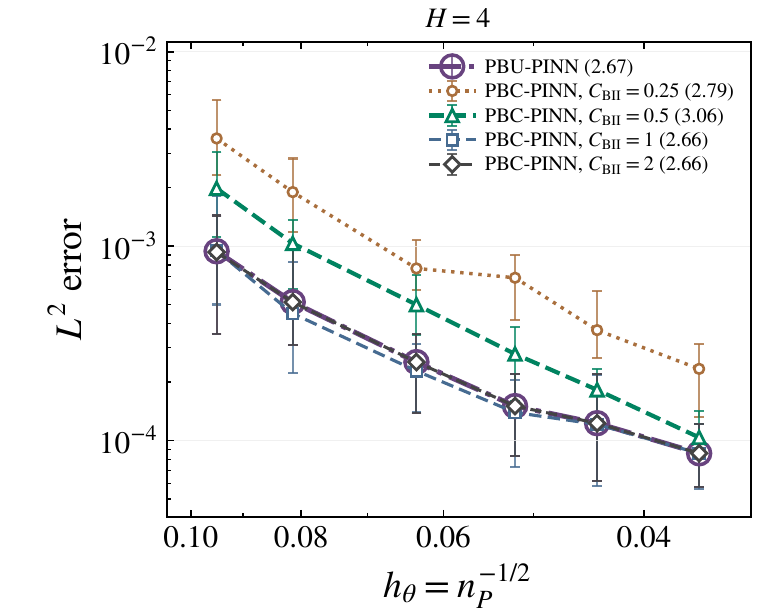}
\end{subfigure}\hspace{0.04\linewidth}%
\begin{subfigure}[t]{0.45\linewidth}
\centering
\includegraphics[width=\linewidth]{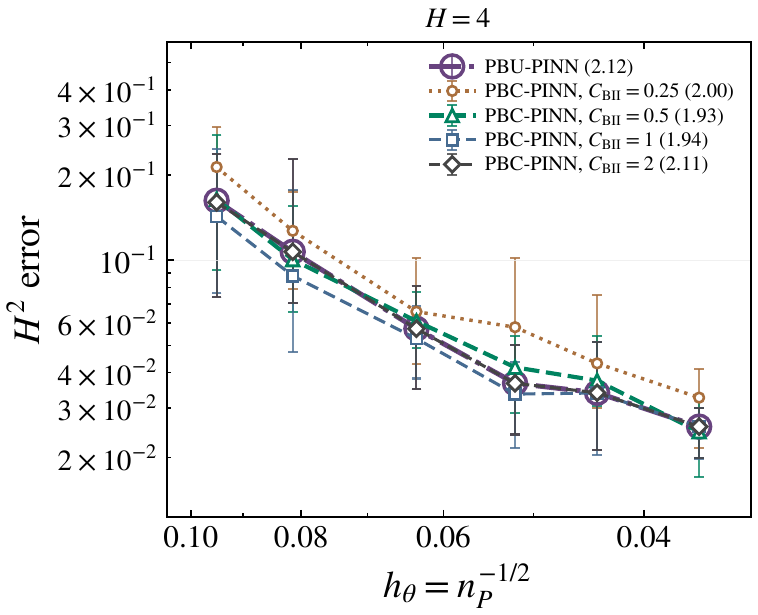}
\end{subfigure}
\caption{Smooth Poisson problem: $L^2$ errors (left) and
$H^2$ errors (right) versus the parameter-based design scale
$h_\theta=n_P^{-1/2}$.
Rows show $H=1,2,4$ hidden layers from top to bottom. PBU-PINN is compared
with PBC-PINN at $C_{\mathrm{BII}}=0.25,0.5,1,2$, with $h_X/h_\theta\approx0.5$
throughout refinement. Curves show ten-seed means with min--max error
bars, and the legend reports fitted convergence orders.}
\label{fig:depth124-cinv-s2}
\vspace{-15pt}
\end{figure}

We next study network refinement on oversampled collocation sets
with a larger training budget. Figure~\ref{fig:depth124-cinv-s2} shows that,
for depth $H=2,4$, PBU-PINN generally achieves smaller errors than PBC-PINN at
$C_{\mathrm{BII}}=0.5$, whereas the latter has higher observed $L^2$ orders
and broadly comparable $H^2$ orders.
Larger BII thresholds generally yield smaller errors, closer to those of
PBU-PINN, but do not consistently improve the observed orders.
Width refinement gives clearer $H^2$ decay with two and four
hidden layers than with one, indicating the importance of sufficient
network depth.

\begin{figure}[t]
\centering
\begin{subfigure}[t]{0.42\linewidth}
\centering
\includegraphics[width=\linewidth]{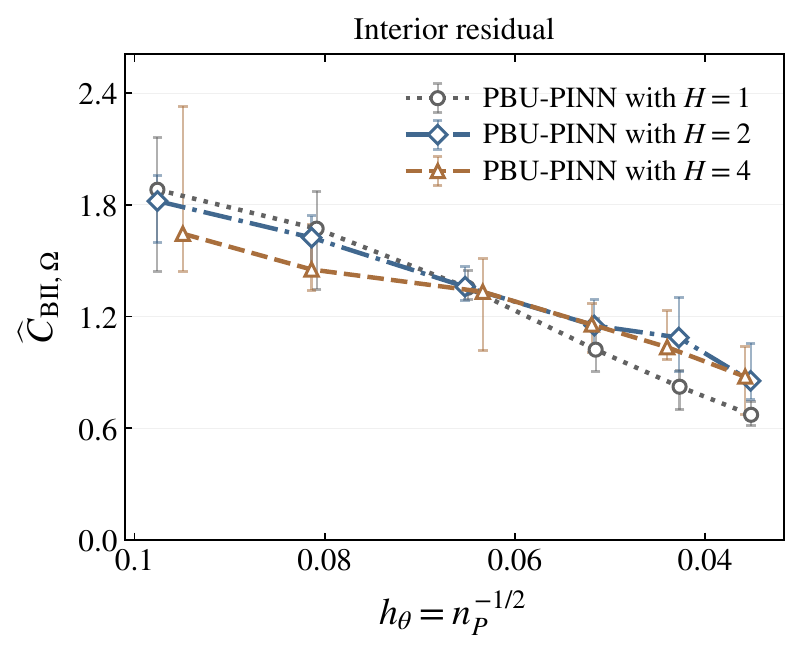}
\end{subfigure}\hspace{0.04\linewidth}%
\begin{subfigure}[t]{0.42\linewidth}
\centering
\includegraphics[width=\linewidth]{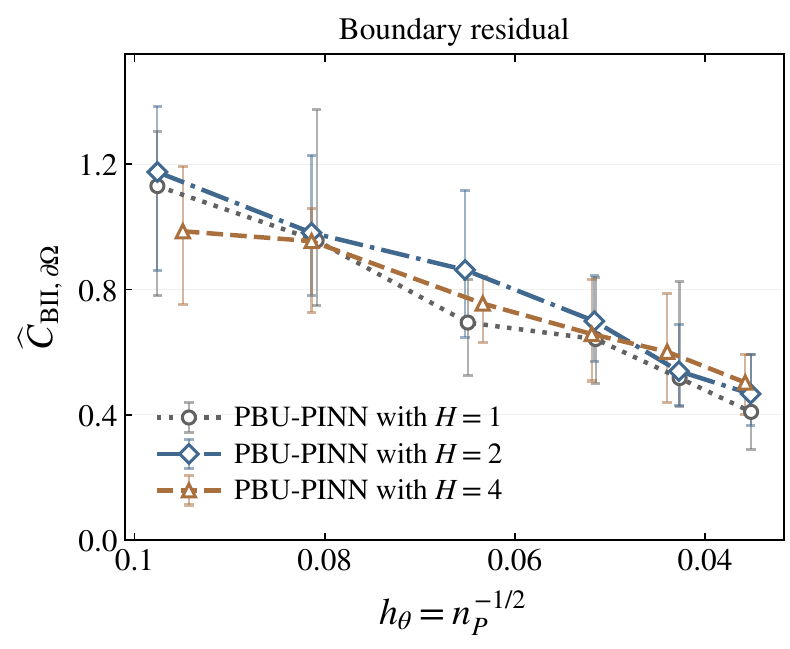}
\end{subfigure}
\caption{Smooth Poisson problem: scaled empirical BII ratios for
the trained PBU-PINNs in Fig.~\ref{fig:depth124-cinv-s2} with $H=1,2,4$ hidden
layers. The left and right panels correspond to the interior and boundary
residuals, respectively. The plots show $\widehat C_{\mathrm{BII},D}=h_\theta R_D$
versus $h_\theta$, where $R_D$ is the stabilized ratio of the empirical $H^2$ and $H^1$
norms of the residual. Curves show ten-seed means with min--max error bars.}
\label{fig:pb-bii-refinement}
\vspace{-15pt}
\end{figure}

For PBU-PINN and $D\in\{\Omega,\partial\Omega\}$, define the scaled empirical BII ratio by
\begin{equation}\label{eq:pb-empirical-bii-requirement}
\widehat C_{\mathrm{BII},D}(h_\theta)
:=h_\theta R_D(v_\theta^{\mathrm{PBU}_\gamma}),
\end{equation}
where $v_\theta^{\mathrm{PBU}_\gamma}$ denotes the trained PBU network and $R_D$ is the stabilized, quadrature-based ratio in
\eqref{eq:implemented-bii-fraction} with $\beta_D=2$.
In Figure~\ref{fig:pb-bii-refinement}, these quantities generally
decrease under refinement over the tested widths,
consistent with the role of bounded BII constants in
Theorems~\ref{thm:discrete} and~\ref{thm:pb-conditional-convergence}.
At $C_{\mathrm{BII}}=2$, the PBC-PINN curves nearly overlap those of PBU-PINN.
These quadrature-based diagnostics do not establish that the continuous BII holds uniformly as $h_\theta\to0$.

% \newpage

\subsection{\texorpdfstring{Example 2. Convergence for finite-regularity solutions}{Example 2. Convergence for finite-regularity solutions}}
\label{subsec:finite-regularity-poisson}

The convergence estimates depend on the solution regularity. We therefore
examine how reduced regularity affects the observed convergence under
network refinement.
On $\Omega=(0,1)^2$, we consider
$
-\Delta u=f_{m_\star}
$
with
boundary data $g_{m_\star}=u_{m_\star}^*|_{\partial\Omega}$ and exact solutions
\begin{equation}
\label{eq:finite-regularity-exact-solution}
u_{m_\star}^*(x,y)
=A_{m_\star}(1-x^2)^2(1-y^2)^2 r^{m_\star-1}\cos\!\big((m_\star-1)\theta\big),
\end{equation}
where $r=\sqrt{x^2+y^2}$ and $\theta=\arctan(y/x)\in(0,\pi/2)$.
We take $m_\star\in\{2.6,3.6,4.6\}$, set $f_{m_\star}=-\Delta u_{m_\star}^*$, and choose $A_{m_\star}>0$ so that
$\|f_{m_\star}\|_{L^2(\Omega)}=\pi^2$.
Here $m_\star$ is the critical Sobolev regularity index:
\begin{align*}
u_{m_\star}^*\in H^{m_\star-\varepsilon}(\Omega)
\quad\text{for every }\varepsilon>0,
\qquad\text{but}\qquad
u_{m_\star}^*\notin H^{m_\star}(\Omega).
\end{align*}
We sample the interior more densely than the boundary.

\begin{figure}[t]
\centering
\hspace{-0.02\linewidth}
\begin{subfigure}[t]{0.33\linewidth}
\centering
\includegraphics[width=\linewidth]{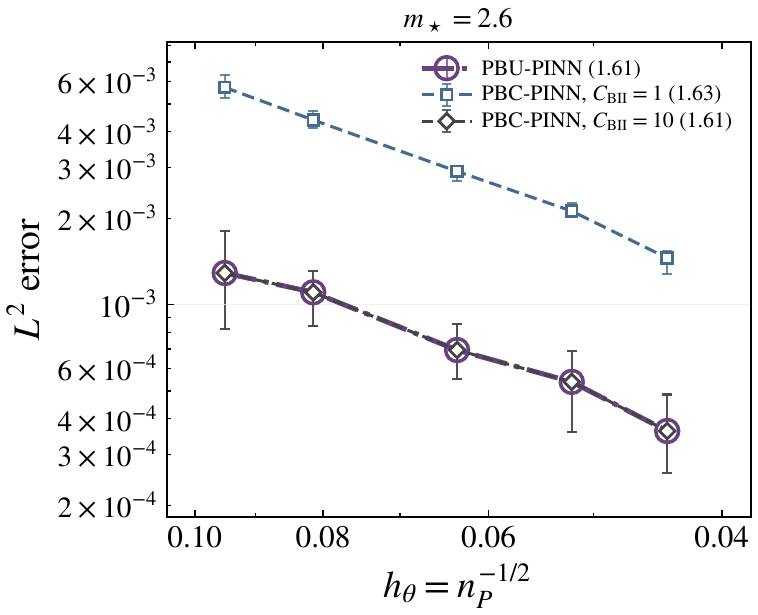}
\end{subfigure}\hspace{-0.001\linewidth}%
\begin{subfigure}[t]{0.33\linewidth}
\centering
\includegraphics[width=\linewidth]{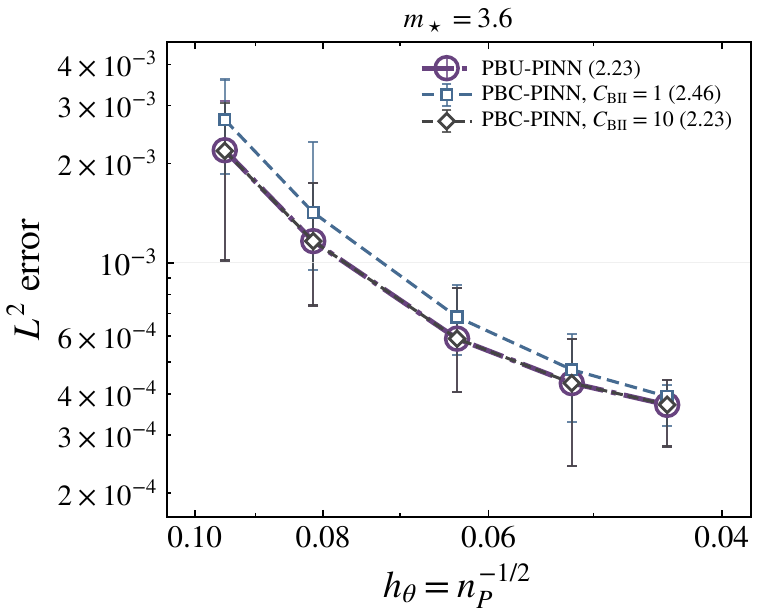}
\end{subfigure}\hspace{-0.001\linewidth}%
\begin{subfigure}[t]{0.33\linewidth}
\centering
\includegraphics[width=\linewidth]{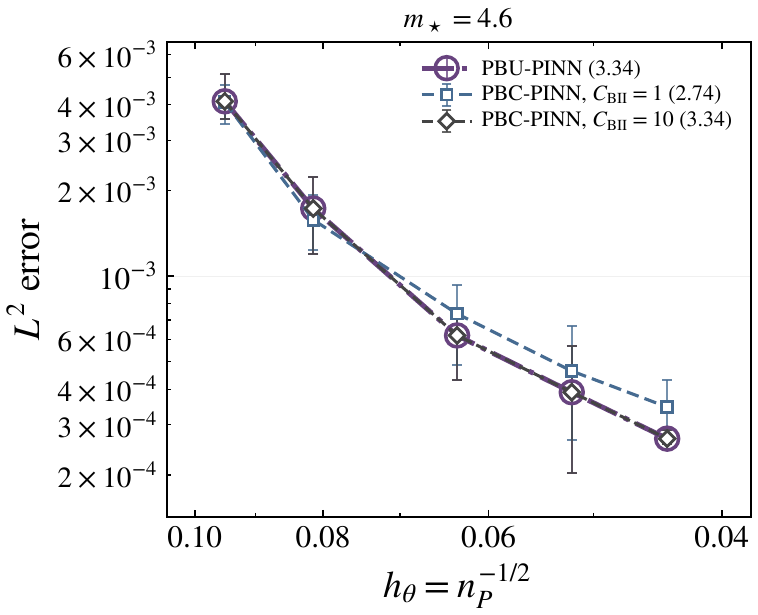}
\end{subfigure}
\par\smallskip
\hspace{-0.02\linewidth}
\begin{subfigure}[t]{0.33\linewidth}
\centering
\includegraphics[width=\linewidth]{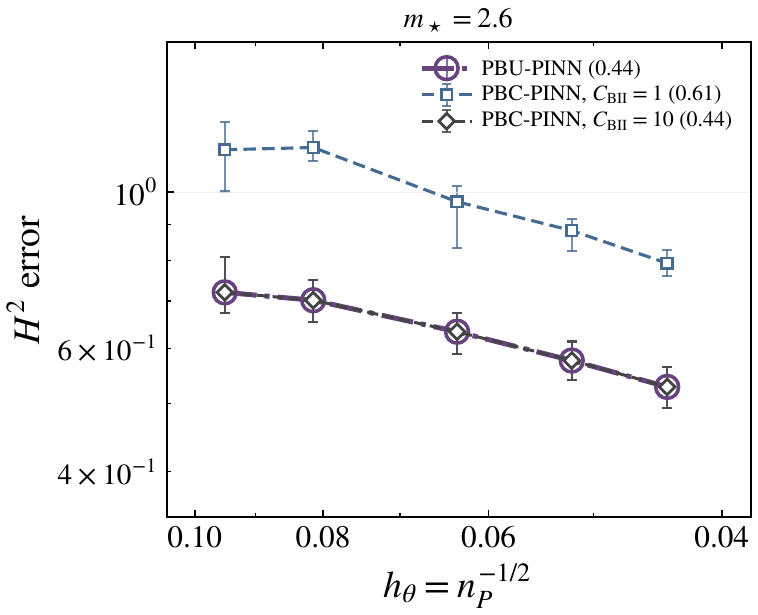}
\end{subfigure}\hspace{-0.001\linewidth}%
\begin{subfigure}[t]{0.33\linewidth}
\centering
\includegraphics[width=\linewidth]{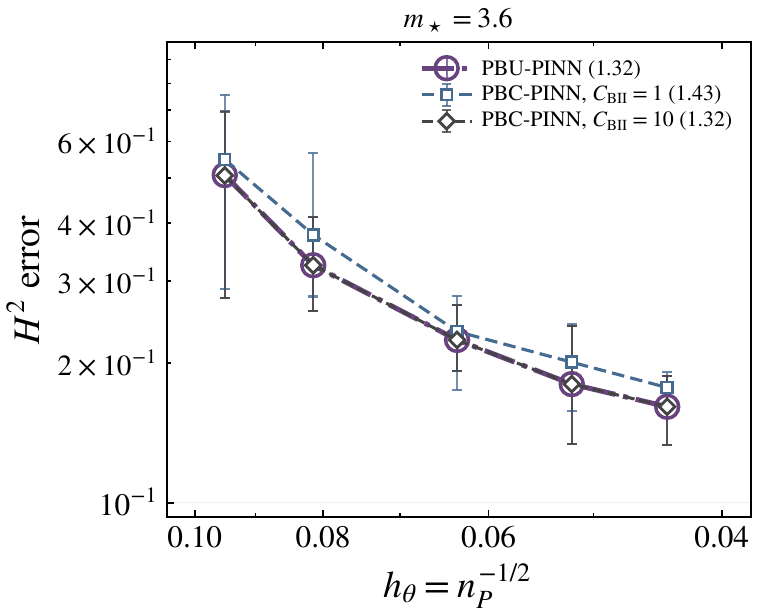}
\end{subfigure}\hspace{-0.001\linewidth}%
\begin{subfigure}[t]{0.33\linewidth}
\centering
\includegraphics[width=\linewidth]{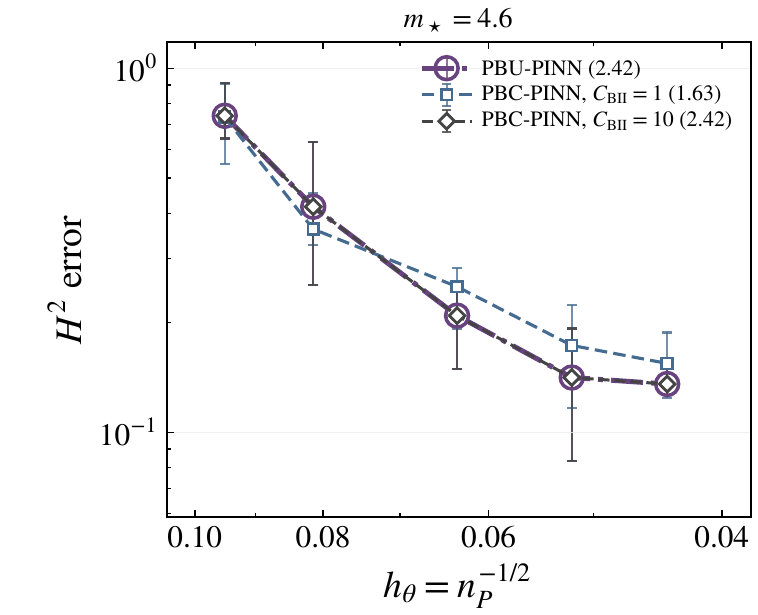}
\end{subfigure}
\caption{Poisson problems with finite solution regularity:
$L^2$ errors (top) and $H^2$ errors (bottom) versus the parameter-based design
scale $h_\theta=n_P^{-1/2}$ for networks with four hidden layers.
The columns show solution regularities $m_\star=2.6,3.6,4.6$ from left to right. PBU-PINN is compared with
PBC-PINN at $C_{\mathrm{BII}}=1,10$, with $h_X/h_\theta\approx0.1$ throughout
refinement. Curves show five-seed means with min--max error bars, and the
legend reports fitted convergence orders.}
\label{fig:finite-regularity-network-refinement}
% \vspace{-15pt}
% \end{figure}

% \begin{figure}[t]
\centering
\begin{subfigure}[t]{0.42\linewidth}
\centering
\includegraphics[width=\linewidth]{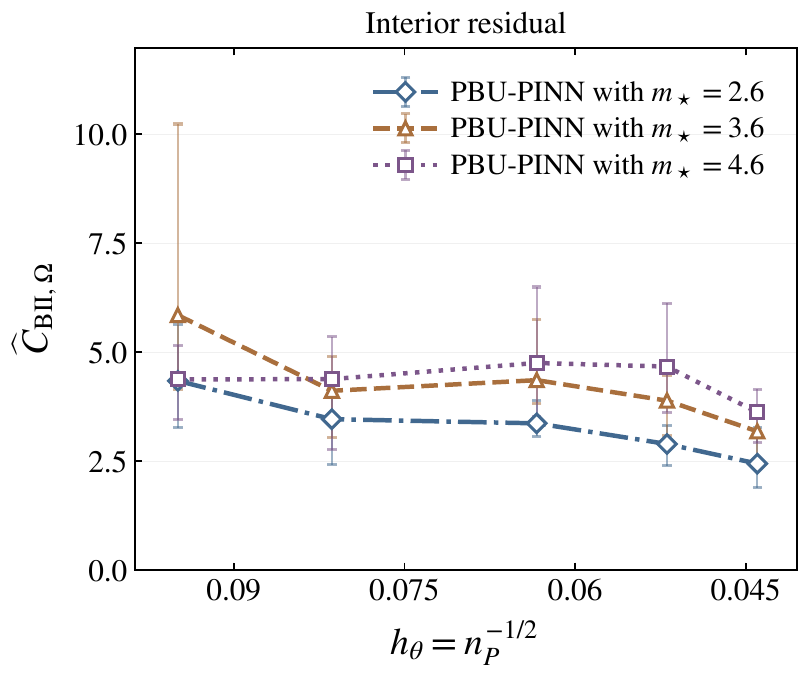}
\end{subfigure}\hspace{0.04\linewidth}%
\begin{subfigure}[t]{0.42\linewidth}
\centering
\includegraphics[width=\linewidth]{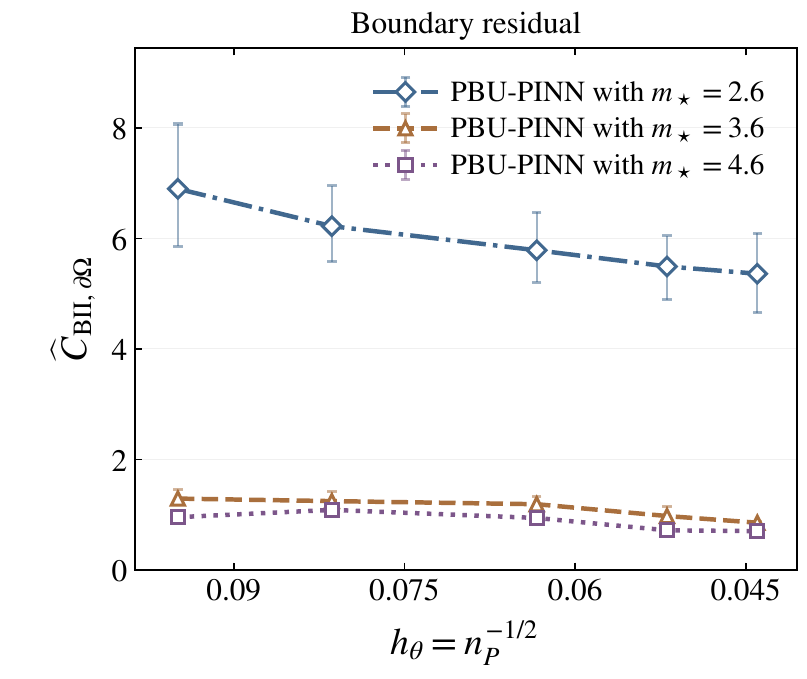}
\end{subfigure}
\caption{Poisson problems with finite solution regularity:
scaled empirical BII ratios for PBU-PINN in Fig.~\ref{fig:finite-regularity-network-refinement} with four hidden layers at solution
regularities $m_\star=2.6,3.6,4.6$. The left and right
panels show $\widehat C_{\mathrm{BII},D}=h_\theta R_D$ versus $h_\theta$ for the interior and
boundary residuals, respectively, where $R_D$ is the stabilized ratio of the
empirical $H^2$ and $H^1$ norms of the residual. Curves show five-seed means with min--max error bars.}
\label{fig:finite-pb-bii-refinement}
\vspace{-15pt}
\end{figure}

Figure~\ref{fig:finite-regularity-network-refinement} shows error decay
in all three cases. The observed $H^2$ convergence orders increase with
solution regularity. The cases with $m_\star=2.6$ and $3.6$ fall outside the regularity assumptions
of Theorem~\ref{thm:discrete}, so the comparison is qualitative.
At $C_{\mathrm{BII}}=1$, PBC-PINN generally has larger errors than PBU-PINN,
while its observed orders relative to PBU-PINN vary with regularity.
The larger error gap at $m_\star=2.6$ may reflect stronger competition between
fitting the corner-singular boundary data and controlling the
boundary-residual BII ratio.
At $C_{\mathrm{BII}}=1$, the final PBC-PINN networks exceeded the
empirical BII threshold.
Figure~\ref{fig:finite-pb-bii-refinement} shows no rapid growth
in the ratios \eqref{eq:pb-empirical-bii-requirement} over the tested widths.
At $C_{\mathrm{BII}}=10$, the PBC curves nearly coincide with PBU-PINN.

\subsection{\texorpdfstring{Example 3. Error decay under network refinement for a semilinear problem}{Example 3. Error decay under network refinement for a semilinear problem}}
The analysis above concerns linear elliptic problems. As a numerical
extension, we consider the semilinear problem
\begin{align}
-\Delta u+10u^3=f
\end{align} 
on $\Omega=(0,1)^2$ with homogeneous Dirichlet data.
The exact solution is $u^*(x,y)=\sin(\pi x)\sin(\pi y)$, and the prescribed
forcing is $f(x,y)=2\pi^2u^*(x,y)+10u^*(x,y)^3$.
In this example, the full residuals are
\begin{align}
r_\Omega(v_\theta)=-\Delta v_\theta+10v_\theta^3-f  \quad \text{and} \quad
r_{\partial\Omega}(v_\theta)=v_\theta,
\end{align}
PBC-PINN applies the penalties
$\mathcal B_D^\delta(v_\theta)$, $D\in\{\Omega,\partial\Omega\}$, to these residuals
in \eqref{eq:total-training-objective}.
The architectures and sampling rules follow Figure~\ref{fig:depth124-cinv-s2},
with paired runs sharing initial parameters and collocation points.
Training uses 2000 Adam steps followed by 80 L-BFGS blocks.

\begin{figure}[t]
\centering
\begin{subfigure}[t]{0.45\linewidth}
\centering
\includegraphics[width=\linewidth]{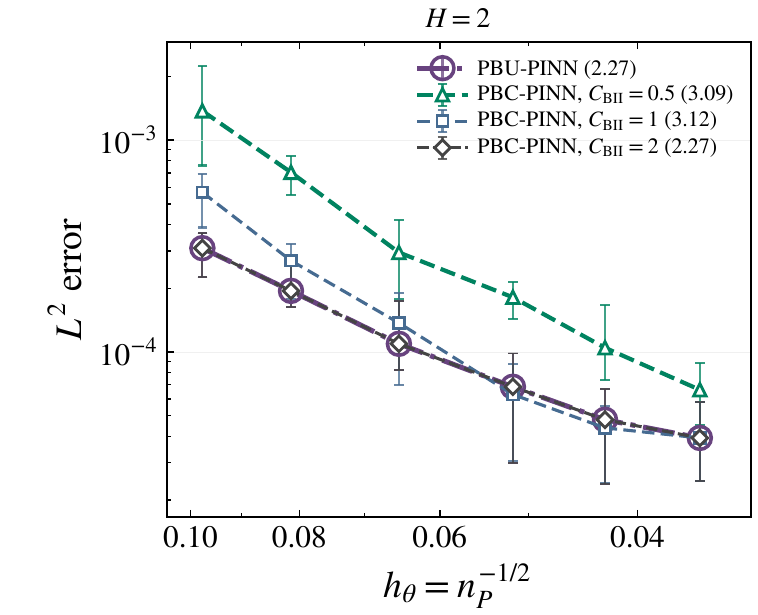}
\end{subfigure}\hspace{0.04\linewidth}%
\begin{subfigure}[t]{0.45\linewidth}
\centering
\includegraphics[width=\linewidth]{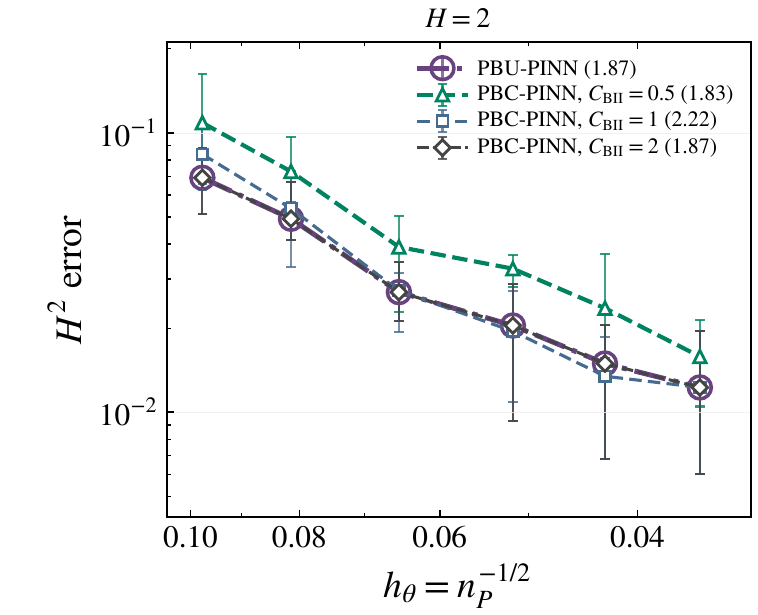}
\end{subfigure}
\par\smallskip
\begin{subfigure}[t]{0.45\linewidth}
\centering
\includegraphics[width=\linewidth]{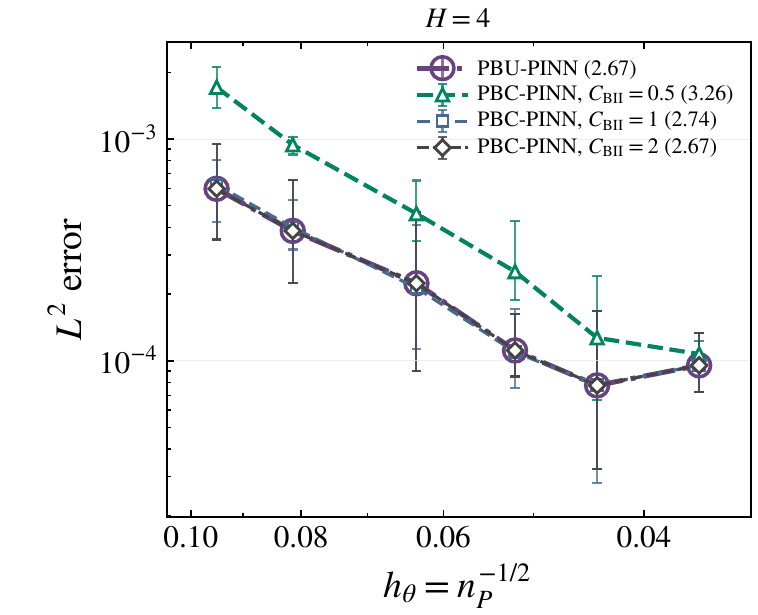}
\end{subfigure}\hspace{0.04\linewidth}%
\begin{subfigure}[t]{0.45\linewidth}
\centering
\includegraphics[width=\linewidth]{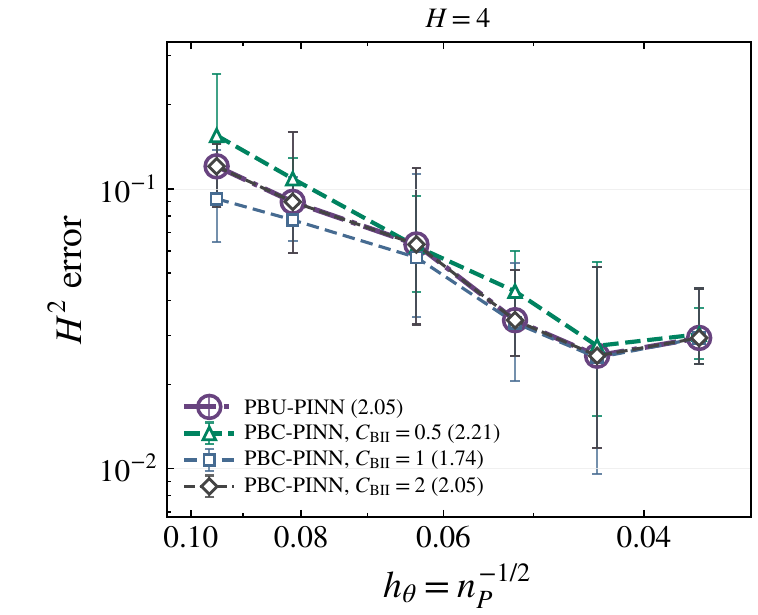}
\end{subfigure}
\caption{Semilinear problem: $L^2$ errors (left) and
$H^2$ errors (right) versus the design scale
$h_\theta=n_P^{-1/2}$.
Rows show $H=2,4$ hidden layers from top to bottom. PBU-PINN is compared
with PBC-PINN at $C_{\mathrm{BII}}=0.5,1,2$, with $h_X/h_\theta\approx0.5$
throughout refinement. }
\label{fig:nonlinear-network-refinement}
% \vspace{-15pt}
% \end{figure}

% \begin{figure}[t]
\centering
\begin{subfigure}[t]{0.41\linewidth}
\centering
\includegraphics[width=\linewidth]{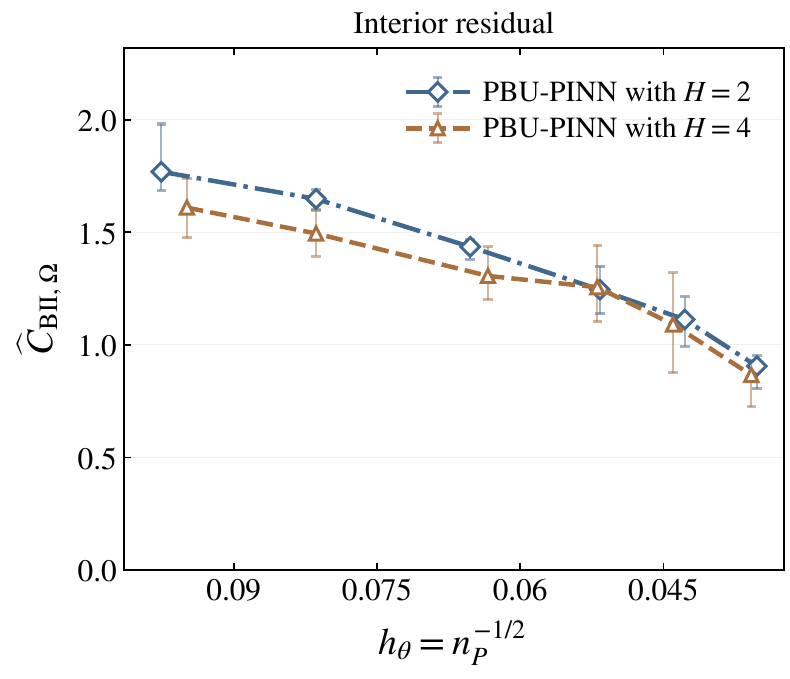}
\end{subfigure}\hspace{0.04\linewidth}%
\begin{subfigure}[t]{0.41\linewidth}
\centering
\includegraphics[width=\linewidth]{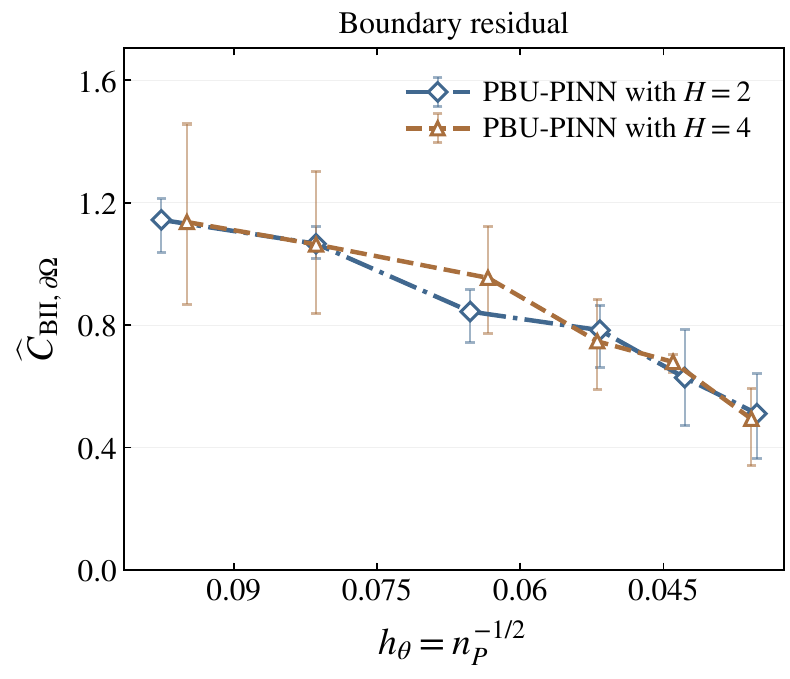}
\end{subfigure}
\caption{Semilinear problem: scaled empirical BII ratios for
PBU-PINN with $H=2,4$ hidden layers in Fig.~\ref{fig:nonlinear-network-refinement}. The left and right
panels show $\widehat C_{\mathrm{BII},D}=h_\theta R_D$ versus $h_\theta$ for the interior and
boundary residuals, respectively, where $R_D$ is the stabilized ratio of the
empirical $H^2$ and $H^1$ norms of the residual.}
\label{fig:nonlinear-pb-bii-refinement}
\vspace{-15pt}
\end{figure}

Figure~\ref{fig:nonlinear-network-refinement} shows that, at the same network size,
PBC-PINN with $C_{\mathrm{BII}}=0.5$ generally has larger $L^2$ errors than PBU-PINN
but exhibits faster error decay under refinement.
The two methods have comparable $H^2$ decay.
Increasing depth from two to four hidden layers does not consistently improve
the observed $H^2$ orders across the tested BII thresholds.
In Figure~\ref{fig:nonlinear-pb-bii-refinement}, the ratios
\eqref{eq:pb-empirical-bii-requirement} generally decrease under refinement
over the tested widths.
At $C_{\mathrm{BII}}=2$, no penalty activation is recorded during paired
training, and the resulting networks coincide with PBU-PINN, explaining
the overlapping convergence curves.
Future work will focus on establishing stability and consistency for the corresponding constrained neural subclasses.

\FloatBarrier
\section{Conclusion}
The stability and consistency developed for kernel-based collocation continue to provide useful guidance for neural methods based on residual minimization. For linear elliptic Dirichlet problems, we have extended this framework to neural trial classes by identifying admissible subclasses on which sampled residuals control the continuous solution error. Inverse estimates for the residuals provide the bounds needed for discrete stability. The analysis also determines the scale dependence of the physics-balanced boundary weight.
For PBC residual minimization, admissibility follows from restricting the trial class to networks that satisfy the BII constraints.
At any fixed network size, every feasible network is admissible
by definition. Under the sampling condition, stability holds.
This alone does not guarantee a nonempty class, an accurate
competitor, or existence of a minimizer.
For any fixed choice of positive BII thresholds, the $\tanh$-MLP construction supplies accurate admissible comparison networks for all sufficiently small $h_\theta$. The stability restriction therefore preserves approximation under refinement, establishing asymptotic consistency for PBC.
For PBU residual minimization, eventual admissibility is conditional. If the scaled residual BII ratios of the minimizers have finite $\limsup$, thresholds chosen above these limits make the minimizers eventually admissible. With sufficiently dense sampling, stability then applies, and minimality requires only a comparison network from the full trial class.
In extending the present stability analysis to shallow $\mathrm{ReLU}^k$
networks, recent work on Bernstein inequalities~\cite{LiuMaoXu2026DiscreteLS,
LiuMaoXu2026Geometry} offers a natural starting point for establishing
the corresponding residual inverse estimates.

The numerical experiments illustrate how network refinement and the BII thresholds affect the error in the linear problems. The nonlinear example exhibits similar error decay and similar trends in the scaled empirical BII ratios.  A corresponding theory for nonlinear problems remains to be established, and the observed similarities provide guidance for further work on stability and consistency within the associated admissible neural classes.

\FloatBarrier
\bibliographystyle{siamplain-nodoi}
\bibliography{RBF}

\end{document}